%% file: 572.tex
\ifx\shlhetal\undefinedcontrolsequence\let\shlhetal\relax\fi

\input amstex
\input mathdefs
\localtags
\NoBlackBoxes
\documentstyle {amsppt}
\topmatter
\title {Colouring and non-productivity of $\aleph_2$-c.c. \\
Sh572} \endtitle
\author {Saharon Shelah \thanks{\null\newline
\S1,\S2 Done 9/94;\S3 written 10/94;\S4 written 1/95. \null\newline
I would like to thank Alice Leonhardt for the beautiful typing and
Zoran Spasojevic for helping to proof read it. \null\newline
Latest Revision 96/Jan/19} \endthanks} \endauthor
\affil {Institute of Mathematics \\
The Hebrew University \\
Jerusalem, Israel
\medskip
Rutgers University \\
Department of Mathematics \\
New Brunswick, NJ  USA} \endaffil
\abstract  We prove that colouring of pairs from $\aleph_2$ with 
strong properties exists.
The easiest to state (and quite a well known problem) it solves: there are
two topological spaces with cellularity $\aleph_1$ whose product has
cellularity $\aleph_2$; equivalently we can speak on cellularity of Boolean
algebras or on Boolean algebras satisfying the $\aleph_2$-c.c. whose product
fails the $\aleph_2$-c.c.  We also deal more with guessing of clubs.
\endabstract

\endtopmatter
\document  

\def\renewcommand{\newcommand}	       
\edef\cite{\the\catcode`@}%
\catcode`@ = 11
\let\@oldatcatcode = \cite
\chardef\@letter = 11
\chardef\@other = 12
%
%
%
%
\def\@innerdef#1#2{\edef#1{\expandafter\noexpand\csname #2\endcsname}}%
%
%
\@innerdef\@innernewcount{newcount}%
\@innerdef\@innernewdimen{newdimen}%
\@innerdef\@innernewif{newif}%
\@innerdef\@innernewwrite{newwrite}%
%
%
%
\def\@gobble#1{}%
%
%
%
\ifx\inputlineno\@undefined
   \let\@linenumber = \empty 
\else
   \def\@linenumber{\the\inputlineno:\space}%
\fi
%
%
%
\def\@futurenonspacelet#1{\def\cs{#1}%
   \afterassignment\@stepone\let\@nexttoken=
}%
\begingroup 
\def\\{\global\let\@stoken= }%
\\ 
\endgroup
\def\@stepone{\expandafter\futurelet\cs\@steptwo}%
\def\@steptwo{\expandafter\ifx\cs\@stoken\let\@@next=\@stepthree
   \else\let\@@next=\@nexttoken\fi \@@next}%
\def\@stepthree{\afterassignment\@stepone\let\@@next= }%
%
%
%
\def\@getoptionalarg#1{%
   \let\@optionaltemp = #1%
   \let\@optionalnext = \relax
   \@futurenonspacelet\@optionalnext\@bracketcheck
}%
%
%
\def\@bracketcheck{%
   \ifx [\@optionalnext
      \expandafter\@@getoptionalarg
   \else
      \let\@optionalarg = \empty
      \expandafter\@optionaltemp
   \fi
}%
\def\@@getoptionalarg[#1]{%
   \def\@optionalarg{#1}%
   \@optionaltemp
}%
%
%
%
\def\@nnil{\@nil}%
\def\@fornoop#1\@@#2#3{}%
\def\@for#1:=#2\do#3{%
   \edef\@fortmp{#2}%
   \ifx\@fortmp\empty \else
      \expandafter\@forloop#2,\@nil,\@nil\@@#1{#3}%
   \fi
}%
\def\@forloop#1,#2,#3\@@#4#5{\def#4{#1}\ifx #4\@nnil \else
       #5\def#4{#2}\ifx #4\@nnil \else#5\@iforloop #3\@@#4{#5}\fi\fi
}%
\def\@iforloop#1,#2\@@#3#4{\def#3{#1}\ifx #3\@nnil
       \let\@nextwhile=\@fornoop \else
      #4\relax\let\@nextwhile=\@iforloop\fi\@nextwhile#2\@@#3{#4}%
}%
%
%
%
\@innernewif\if@fileexists
\def\@testfileexistence{\@getoptionalarg\@finishtestfileexistence}%
\def\@finishtestfileexistence#1{%
   \begingroup
      \def\extension{#1}%
      \immediate\openin0 =
         \ifx\@optionalarg\empty\jobname\else\@optionalarg\fi
         \ifx\extension\empty \else .#1\fi
         \space
      \ifeof 0
         \global\@fileexistsfalse
      \else
         \global\@fileexiststrue
      \fi
      \immediate\closein0
   \endgroup
}%
%
%
%
%
\def\bibliographystyle#1{%
   \@readauxfile
   \@writeaux{\string\bibstyle{#1}}%
}%
\let\bibstyle = \@gobble
%
%
\let\bblfilebasename = \jobname
\def\bibliography#1{%
   \@readauxfile
   \@writeaux{\string\bibdata{#1}}%
   \@testfileexistence[\bblfilebasename]{bbl}%
   \if@fileexists
      \nobreak
      \@readbblfile
   \fi
}%
\let\bibdata = \@gobble
%
%
\def\nocite#1{%
   \@readauxfile
   \@writeaux{\string\citation{#1}}%
}%
\@innernewif\if@notfirstcitation
%
%
\def\cite{\@getoptionalarg\@cite}%
%
%
\def\@cite#1{%
   \let\@citenotetext = \@optionalarg
   \printcitestart
   \nocite{#1}%
   \@notfirstcitationfalse
   \@for \@citation :=#1\do
   {%
      \expandafter\@onecitation\@citation\@@
   }%
   \ifx\empty\@citenotetext\else
      \printcitenote{\@citenotetext}%
   \fi
   \printcitefinish
}%
\def\@onecitation#1\@@{%
   \if@notfirstcitation
      \printbetweencitations
   \fi
   \expandafter \ifx \csname\@citelabel{#1}\endcsname \relax
      \if@citewarning
         \message{\@linenumber Undefined citation `#1'.}%
      \fi
      \expandafter\gdef\csname\@citelabel{#1}\endcsname{%
\strut
\vadjust{\vskip-\dp\strutbox
\vbox to 0pt{\vss\parindent0cm \leftskip=\hsize 
\advance\leftskip3mm
\advance\hsize 4cm\strut\openup-4pt 
\rightskip 0cm plus 1cm minus 0.5cm ?  #1 ?\strut}}
         {\tt
            \escapechar = -1
            \nobreak\hskip0pt
            \expandafter\string\csname#1\endcsname
            \nobreak\hskip0pt
         }%
      }%
   \fi
   \csname\@citelabel{#1}\endcsname
   \@notfirstcitationtrue
}%
%
%
\def\@citelabel#1{b@#1}%
%
%
\def\@citedef#1#2{\expandafter\gdef\csname\@citelabel{#1}\endcsname{#2}}%
%
%
%
\def\@readbblfile{%
   \ifx\@itemnum\@undefined
      \@innernewcount\@itemnum
   \fi
   \begingroup
      \def\begin##1##2{%
         \setbox0 = \hbox{\biblabelcontents{##2}}%
         \biblabelwidth = \wd0
      }%
      \def\end##1{}
      %
      %
      \@itemnum = 0
      \def\bibitem{\@getoptionalarg\@bibitem}%
      \def\@bibitem{%
         \ifx\@optionalarg\empty
            \expandafter\@numberedbibitem
         \else
            \expandafter\@alphabibitem
         \fi
      }%
      \def\@alphabibitem##1{%
         \expandafter \xdef\csname\@citelabel{##1}\endcsname {\@optionalarg}%
         \ifx\biblabelprecontents\@undefined
            \let\biblabelprecontents = \relax
         \fi
         \ifx\biblabelpostcontents\@undefined
            \let\biblabelpostcontents = \hss
         \fi
         \@finishbibitem{##1}%
      }%
      \def\@numberedbibitem##1{%
         \advance\@itemnum by 1
         \expandafter \xdef\csname\@citelabel{##1}\endcsname{\number\@itemnum}%
         \ifx\biblabelprecontents\@undefined
            \let\biblabelprecontents = \hss
         \fi
         \ifx\biblabelpostcontents\@undefined
            \let\biblabelpostcontents = \relax
         \fi
         \@finishbibitem{##1}%
      }%
      \def\@finishbibitem##1{%
         \biblabelprint{\csname\@citelabel{##1}\endcsname}%
         \@writeaux{\string\@citedef{##1}{\csname\@citelabel{##1}\endcsname}}%
         \ignorespaces
      }%
      %
      %
      \let\em = \bblem
      \let\newblock = \bblnewblock
      \let\sc = \bblsc
      \frenchspacing
      \clubpenalty = 4000 \widowpenalty = 4000
      \tolerance = 10000 \hfuzz = .5pt
      \everypar = {\hangindent = \biblabelwidth
                      \advance\hangindent by \biblabelextraspace}%
      \bblrm
      \parskip = 1.5ex plus .5ex minus .5ex
      \biblabelextraspace = .5em
      \bblhook
      \input \bblfilebasename.bbl
   \endgroup
}%
%
%
\@innernewdimen\biblabelwidth
\@innernewdimen\biblabelextraspace
%
%
%
\def\biblabelprint#1{%
   \noindent
   \hbox to \biblabelwidth{%
      \biblabelprecontents
      \biblabelcontents{#1}%
      \biblabelpostcontents
   }%
   \kern\biblabelextraspace
}%
%
%
%
\def\biblabelcontents#1{{\bblrm [#1]}}%
%
%
\def\bblrm{\rm}%
%
%
\def\bblem{\it}%
%
%
\def\bblsc{\ifx\@scfont\@undefined
              \font\@scfont = cmcsc10
           \fi
           \@scfont
}%
%
%
\def\bblnewblock{\hskip .11em plus .33em minus .07em }%
%
%
\let\bblhook = \empty
%
%
%
\def\printcitestart{[}
\def\printcitefinish{]}
\def\printbetweencitations{, }
\def\printcitenote#1{, #1}
%
%
%
\let\citation = \@gobble
%
%
%
\@innernewcount\@numparams
%
%
\def\newcommand#1{%
   \def\@commandname{#1}%
   \@getoptionalarg\@continuenewcommand
}%
%
%
\def\@continuenewcommand{%
   \@numparams = \ifx\@optionalarg\empty 0\else\@optionalarg \fi \relax
   \@newcommand
}%
%
%
\def\@newcommand#1{%
   \def\@startdef{\expandafter\edef\@commandname}%
   \ifnum\@numparams=0
      \let\@paramdef = \empty
   \else
      \ifnum\@numparams>9
         \errmessage{\the\@numparams\space is too many parameters}%
      \else
         \ifnum\@numparams<0
            \errmessage{\the\@numparams\space is too few parameters}%
         \else
            \edef\@paramdef{%
               \ifcase\@numparams
                  \empty  No arguments.
               \or ####1%
               \or ####1####2%
               \or ####1####2####3%
               \or ####1####2####3####4%
               \or ####1####2####3####4####5%
               \or ####1####2####3####4####5####6%
               \or ####1####2####3####4####5####6####7%
               \or ####1####2####3####4####5####6####7####8%
               \or ####1####2####3####4####5####6####7####8####9%
               \fi
            }%
         \fi
      \fi
   \fi
   \expandafter\@startdef\@paramdef{#1}%
}%
%
%
%
%
\def\@readauxfile{%
   \if@auxfiledone \else 
      \global\@auxfiledonetrue
      \@testfileexistence{aux}%
      \if@fileexists
         \begingroup
            \endlinechar = -1
            \catcode`@ = 11
            \input \jobname.aux
         \endgroup
      \else
         \message{\@undefinedmessage}%
         \global\@citewarningfalse
      \fi
      \immediate\openout\@auxfile = \jobname.aux
   \fi
}%
%
%
\newif\if@auxfiledone
\ifx\noauxfile\@undefined \else \@auxfiledonetrue\fi
%
%
%
%
\@innernewwrite\@auxfile
\def\@writeaux#1{\ifx\noauxfile\@undefined \write\@auxfile{#1}\fi}%
%
%
%
\ifx\@undefinedmessage\@undefined
   \def\@undefinedmessage{No .aux file; I won't give you warnings about
                          undefined citations.}%
\fi
%
%
\@innernewif\if@citewarning
\ifx\noauxfile\@undefined \@citewarningtrue\fi
%
%
%
\catcode`@ = \@oldatcatcode
\newpage
\sectno=-1
\head {Annotated content} \endhead
\resetall
\medskip
\roster
\item "{\S1}"  Retry at $\aleph_2$-c.c. not productive  \newline
\smallskip
\noindent
[We prove $Pr_1(\aleph_1,\aleph_2,\aleph_2,\aleph_0)$ which is a much 
stronger result].
\bigskip
\item "{\S2}"  The implicit properties \newline
\smallskip
\noindent
[We define a property implicit in \S1, note what the proof in \S1 gives,
and look at related implication for successor of singular non-strong limit and
show that $Pr_1$ implies $Pr_6$].
\bigskip
\item "{\S3}"  Guessing clubs revisited  \newline
\smallskip
\noindent
[We improve some results mainly from \cite{Sh:413}, giving complete proofs.
We show that for $\mu$ regular uncountable and $\chi < \mu$ we can find
\newline 
$\langle C_\delta:\delta < \mu^+,\text{cf}(\delta) = \mu \rangle$ 
and functions $h_\delta$, from $C_\delta$ onto $\chi$, such that for 
every club $E$ of $\mu^+$ 
for stationarily many $\delta < \mu^+$ we have: $\text{cf}(\delta) = \mu$ 
and for every $\gamma < \chi$ for arbitrarily large $\alpha \in \text{ nacc}
(C_\delta)$ we
have $\alpha \in E,h_\delta(\alpha) = \gamma$.  Also if $C_\delta =
\{ \alpha_{\delta,\varepsilon}:\varepsilon < \mu\}$, $(\alpha_{\delta,
\varepsilon}$ increasing continuous in $\varepsilon$) we can demand
$\{ \varepsilon < \mu:\alpha_{\delta,\varepsilon + 1} \in E$ (and
$\alpha_{\delta,\varepsilon} \in E)\}$ is a stationary subset of $\mu$.  In
fact for each $\gamma < \mu$ the set 
$\{ \varepsilon < \mu:\alpha_{\delta,\varepsilon + 1} \in E,\alpha_{\delta,
\varepsilon} \in E \text{ and } f(\alpha_{\delta,\varepsilon + 1}) = \gamma\}$
is a stationary subset of $\mu$.  We also deal with a parallel to
the last one (without $f$) to successor of singulars and to inaccessibles.]
\bigskip
\item "{\S4}"  More on $Pr_1$ \newline
\smallskip
\noindent
[We prove that Pr$_1(\lambda^{+2},\lambda^{+2},\lambda^{+2},\lambda)$ holds
for regular $\lambda$]. \newline
\medskip
\noindent
On history, references and consequences see \cite[AP1]{Sh:g} and 
\cite[III,\S0]{Sh:g}.
\endroster
\newpage

\head {\S1 Retry at $\aleph_2$-c.c. not productive} \endhead
\resetall
\bigskip

\proclaim{\stag{1.1} Theorem}  $Pr_1(\aleph_2,\aleph_2,\aleph_2,\aleph_0)$.
\endproclaim
\bigskip

\remark{\stag{1.1A} Remark}  1) Is this hard?  
Apostriory it does not look so, but
we have worked hard on it several times without success (worse: produce 
several false proofs).  We thank Juhasz and Soukup for pointing out 
a gap. \newline
2)  Remember that \newline
\underbar{Definition}  $Pr_1(\lambda,\mu,\theta,\sigma)$ means that
there is a symmetric two-place
function $d$ from $\lambda$ to $\theta$ such that: \newline
\underbar{if} $\langle u_\alpha:\alpha < \mu \rangle$ satisfies

$$
u_\alpha \subseteq \lambda,
$$

$$
|u_\alpha| < \sigma,
$$

$$
\alpha < \beta \Rightarrow u_\alpha \cap u_\beta = \emptyset,
$$
\medskip

\noindent
and $\gamma < \theta$ \underbar{then} for some $\alpha < \beta$ we have

$$
\zeta \in u_\alpha \and \xi \in u_\alpha \Rightarrow d(\zeta,\xi) = \gamma.
$$
\medskip

\noindent
3) If we are content with proving that there is a 
colouring with $\aleph_1$ colours, then
we can simplify somewhat: in stage C we let $c(\beta,\alpha) = d_{\text{sq}}
(\rho_{h_1}(\beta,\alpha))$ and this shortens stage D.
\endremark
\bigskip

\demo{Proof}
\enddemo
\bigskip

\noindent
\underbar{Stage A}: First we define a preliminary colouring. \newline
There is a function $d_{sq}:{}^{\omega >}(\omega_1) \rightarrow \omega_1$ 
such that:
\medskip
\roster
\item"{$\bigotimes$}"  \underbar{if} $A \in [\omega_1]^{\aleph_1}$ and
$\langle (\rho_\alpha,\nu_\alpha):\alpha \in A \rangle$ is such that
$\rho_\alpha \in {}^{\omega >}\omega_1,\nu_\alpha \in {}^{\omega >}\omega_1$,
\newline
$\alpha \in \text{ Rang}(\rho_\alpha) \cap \text{ Rang}(\nu_\alpha)$ and
$\gamma < \omega_1$ \underbar{then} for some $\zeta < \xi$ from \newline
$A$ we have: if $\nu',\rho'$ are 
subsequences of $\nu_\zeta,\rho_\xi$ respectively and \newline
$\zeta \in \text{ Rang}(\nu'),\xi \in \text{ Rang}(\rho')$ then
$$
d_{sq}(\nu' \char 94 \rho') = \gamma.
$$
\endroster

\bigskip

\demo{Proof of $\bigotimes$}  
Choose pairwise distinct $\eta_\alpha \in
{}^\omega 2$ for $\alpha < \omega_1$.  Let $d_0:[\omega_1]^2 
\rightarrow \omega_1$ \newline
be such that:
\medskip
\roster
\item "{$(*)$}"  if $n < \omega$ and $\alpha_{\zeta,\ell} < \omega_1$ for
$\zeta < \omega_1,\ell < n$ are pairwise distinct and \newline
$\gamma < \omega_1$
then for some $\zeta < \xi < \omega_1$ we have $\ell < n \Rightarrow \gamma
= d_0(\{\alpha_{\zeta,\ell},\alpha_{\xi,\ell}\})$ \newline
(exists by \cite[see (2.4),p.176]{Sh:261} the $n$ there is $2$).
\endroster
\medskip

\noindent
Define $d_{sq}(\nu)$ for $\nu \in {}^{\omega >}(\omega_1)$ as follows.
If $\ell g(\nu) \le 1$ or $\nu$ is constant then \newline
$d_{sq}(\nu)$ is $0$. Otherwise let \newline
$n(\nu) =: \text{ max}
\{\ell g(\eta_{\nu(\ell)} \cap \eta_{\nu(k)}):\ell < k < \ell g(\nu)$ and
$\nu(\ell) \ne \nu(k)\} < \omega$. \newline
The maximum is on a non-empty set as $\ell g(\nu) \ge 2$ and $\nu$ is not
constant, remember
$\eta_\alpha \in {}^\omega 2$ were pairwise distinct so $\nu(\ell) \ne
\nu(k) \Rightarrow \eta_{\nu(\ell)} \cap \eta_{\nu(k)} \in {}^{\omega >}2$
(is the largest common initial segment of $\eta_{\nu(\ell)},\eta_{\nu(k)}$).
Let $a(\nu) = \{(\ell,k):\ell < k < \ell g(\nu)$ and
$\ell g(\eta_{\nu(\ell)} \cap
\eta_{\nu(k)}) = n(\nu)\}$ so $a(\nu)$ is non-empty and choose the
(lexicographically) minimal pair $(\ell_\nu,k_\nu)$ in it.  Lastly let

$$
d_{sq}(\nu) = d_0(\{ \nu(\ell_\nu),\nu(k_\nu)\}).
$$

\noindent
So $d_{sq}$ is a function with the right domain and range.  Now suppose
we are given $A \in [\omega_1]^{\aleph_1}$,
$\gamma < \omega_1$ and $\rho_\alpha,\nu_\alpha \in
{}^{\omega >}(\omega_1)$ for $\alpha \in A$ such that \newline
$\alpha \in \text{ Rang}(\rho_\alpha) \cap \text{ Rang}(\nu_\alpha)$.
We should find $\alpha < \beta$ from $A$ such that 
$d_{sq}(\nu' \char 94 \rho') = \gamma$ for any subsequences 
$\nu',\rho'$ subsequences of $\nu_\alpha,\rho_\beta$ respectively such
that $\alpha \in \text{ Rang}(\nu')$ and $\beta \in \text{ Rang}(\rho')$.
\medskip

For each $\alpha \in A$ we can find $m_\alpha < \omega$ such that:
\medskip
\roster
\item "{$(*)_0$}"  if $\ell < k < \ell g(\nu_\alpha \char 94 \rho_\alpha)$ and
$(\nu_\alpha \char 94 \rho_\alpha)(\ell) \ne (\nu_\alpha \char 94
\rho_\alpha)(k)$ then \newline
$\eta_{(\nu_\alpha \char 94 \rho_\alpha)(\ell)} \restriction
m_\alpha \ne \eta_{(\nu_\alpha \char 94 \rho_\alpha)(k)} \restriction 
m_\alpha$.
\endroster
\medskip

\noindent
Next we can find $B \in [A]^{\aleph_1}$ such that for all $\alpha \in B$
(the point is that the values do not depend on $\alpha$) we have:
\medskip
\roster
\item "{(a)}"  $\ell g(\nu_\alpha) = m^0,\ell g(\rho_\alpha) = m^1,$
\item "{(b)}"  $a^* = \{(\ell,k):\ell < k < m^0 + m^1$ and 
$(\nu_\alpha \char 94 \rho_\alpha)(\ell) = (\nu_\alpha
\char 94 \rho_\alpha)(k)\},$
\item "{(c)}"  $b^* = \{ \ell < m^0 + m^1:
\alpha = (\nu_\alpha \char 94 \rho_\alpha)(\ell)\}$,
\item "{(d)}"  $m_\alpha = m^2$,
\item "{(e)}"  $\langle \eta_{(\nu_\alpha \char 94 \rho_\alpha)(\ell)}
\restriction m_\alpha:\ell < m^0 + m^1 \rangle = \bar \eta^*$,
\item "{(f)}"  $\langle\text{Rang}(\nu_\alpha \char 94 \rho_\alpha):
\alpha \in B \rangle$ is a $\triangle$-system with heart $w$,
\item "{(g)}"  $u^* = \{ \ell:(\nu_\alpha \char 94 \rho_\alpha)(\ell) \in
w\}$ (so $u^* \ne \{ \ell:\ell < m^0 + m^1\}$ as
$\alpha \in \text{ Rang}(\nu_\alpha \char 94 \rho_\alpha))$,
\item "{(h)}"  $\alpha^*_\ell = (\nu_\alpha \char 94 \rho_\alpha)(\ell)$
for $\ell \in u^*$,
\item "{(i)}"  if $\alpha < \beta \in B$ then $\text{ sup Rang}(\nu_\alpha
\char 94 \rho_\alpha) < \beta$.
\endroster
\medskip

For $\zeta \in B$ let $\bar \beta^\zeta =: \langle (\nu_\zeta \char 94 
\rho_\zeta)(\ell):\ell < m^0 + m^1,\ell \notin u^* \rangle$ and apply 
$(*)$, i.e. the choice of $d_0$.  So for some $\zeta < \xi$ from $B$, we have

$$
\ell < m^0 + m^1 \and \ell \notin u^* \Rightarrow \gamma 
= d_0 \biggl( \{(\nu_\zeta \char 94 
\rho_\zeta)(\ell),(\nu_\xi \char 94 \rho_\xi)(\ell)\} \biggr).
$$ 

We shall prove that $\zeta < \xi$ are as required (in $\otimes$).  So let
$\nu',\rho'$ be subsequences of $\nu_\zeta,\rho_\xi$ (so let
$\nu' = \nu_\zeta \restriction v_1$ and $\rho' = \rho_\xi 
\restriction v_2$) such that $\zeta \in \text{ Rang}(\nu'),\xi \in 
\text{ Rang}(\rho')$ and we have to prove $\gamma = d_{sq}(\nu' 
\char 94 \rho')$.
Let $\tau = \nu' \char 94 \rho'$, so
$\tau = (\nu_\zeta \char 94 \rho_\xi) \restriction (v_1 \cup (m^0 + v_2))$
(in a slight abuse of notation, we look at $\tau$ as a function with domain
$v_1 \cup (m^0 + v_2)$ and also as a member of ${}^{\omega >}(\omega_1)$
where $m + v =: \{ m + \ell:\ell \in v\}$, of course).
By the definition of $d_{sq}$ it is enough to prove the following two things:
\medskip
\roster
%
\item "{$(*)_1$}"  $n(\nu' \char 94 \rho') \ge m^2$ (see clause $(d)$ and
$(*)_0$ above),
\item "{$(*)_2$}"  for every $\ell_1,\ell_2 \in v_1 \cup (m^0 + v_2)$
we have \newline
$\ell g(\eta_{\tau(\ell_1)} \cap
\eta_{\tau(\ell_2)}) \in [m^2,\omega) \Rightarrow
\gamma = d_0(\{\tau(\ell_1),\tau(\ell_2)\})$.
\endroster
\enddemo
\bigskip

\demo{Proof of $(*)_1$}  Let $\ell_1 \in v_1$ and $\ell_2 \in v_2$
be such that $\nu_\zeta(\ell_1) = \zeta$ and $\rho_\xi(\ell_2) = \xi$.

So clearly $\ell_1,m^0 + \ell_2 \in b^*$ (see clause (c)) and \newline
$\eta_{\rho_\xi(\ell_2)} \restriction m^2 = \eta_{\rho_\zeta(\ell_2)} 
\restriction m^2 = \eta_{\nu_\zeta(\ell_1)} \restriction m^2$ (first equality
as $\zeta,\xi \in B$ \newline
and $m_\zeta = m_\xi = m^2$ (see clause (d) and (e)), 
second equality as \newline
$\eta_{\rho_\zeta(\ell_2)} = \eta_{\nu_\zeta(\ell_1)}$ since $\ell_1,m^0 +
\ell_2 \in b^*$ (see clause (c)).  But $\rho_\xi(\ell_2) = \xi \ne 
\zeta = \nu_\zeta(\ell_1)$, hence 
$\eta_{\rho_\xi(\ell_2)} \ne \eta_{\nu_\zeta(\ell_1)}$, so together with the
previous sentence we have

$$
m^2 \le \ell g(\eta_{\nu_\zeta(\ell_1)} \cap \eta_{\rho_\xi(\ell_2)})
= \ell g(\eta_{\tau(\ell_1)} \cap \eta_{\tau(m^0+\ell_2)}) < \omega.
$$
\medskip

\noindent
Hence $n(\tau) \ge m^2$ as required in $(*)_1$.
\enddemo
\bigskip

\demo{Proof of $(*)_2$}  If $\ell_1 < \ell_2$ are from $v_1$, by the choice of
$m^2 = m_\zeta$ it is easy.  Namely, if
$(\ell_1,\ell_2) \in a(\tau)$ then $(\ell_1,\ell_2) \in a(\nu_\zeta)$ and
$\ell g(\eta_{\tau(\ell_1)} \cap \eta_{\tau(\ell_2)}) = 
\ell g(\eta_{\nu_\zeta(\ell_1)} \cap \eta_{\nu_\zeta(\ell_2)}) < m_\zeta 
= m^2$.  
If $\ell_1,\ell_2 \in m^0 + v^2$, by the choice
of $m^2 = m_\xi$ similarly it is easy to show $\ell g(\eta_{\tau(\ell_1)} \cap
\eta_{\tau(\ell_2)}) < m^2$.  So it is enough to prove
\medskip
\roster
\item "{$(*)_3$}"  assume $\ell_1 \in v_1,\ell_2 \in v_2$ and \newline
$\ell g(\eta_{\nu_\zeta(\ell_1)} \cap \eta_{\rho_\xi(\ell_2)}) \in
[m^2,\omega)$ then \newline
$\gamma = d_0(\{ \nu_\zeta(\ell_1),\rho_\xi(\ell_2)\})$.
\endroster
\medskip

\noindent
Now the third assumption in $(*)_3$ means
$\eta_{\nu_\zeta(\ell_1)} \restriction 
m^2 = \eta_{\rho_\xi(\ell_2)} \restriction m^2$ and as 
$\zeta,\xi \in B$ we know that
$\eta_{\rho_\xi(\ell_2)} \restriction m^2 = \eta_{\rho_\zeta(\ell_2)}
\restriction m^2$.  Together we know that $\eta_{\nu_\zeta(\ell_1)} 
\restriction m^2 = 
\eta_{\rho_\zeta(\ell_2)} \restriction m^2$, hence by the choice of
$m_\zeta = m^2$ necessarily
$\eta_{\nu_\zeta(\ell_1)} = \eta_{\rho_\zeta(\ell_2)}$ so that
$\nu_\zeta(\ell_1) = \rho_\zeta(\ell_2)$ and (see clause (b)) also
$\nu_\xi(\ell_1) = \rho_\xi(\ell_2)$.  So

$$
d_0(\{\nu_\zeta(\ell_1),\rho_\xi(\ell_2)\}) = d_0(\{ \nu_\zeta(\ell_1),
\nu_\xi(\ell_1)\}).
$$
\medskip

\noindent
The latter is the required $\gamma$ provided that $\ell_1 \notin u^*$.
Equivalently $\nu_\zeta(\ell_1) \ne \nu_\xi(\ell_1)$ but otherwise also
$\nu_\zeta(\ell_1) = \rho_\xi(\ell_2)$ so $\ell g(\eta_{\nu_\zeta(\ell_1)}
\cap \eta_{\rho_\xi(\ell_2)}) = \omega$, contradicting the assumption of
$(*)_3$ that $\ell g(\eta_{\tau(\ell_1)} \cap \eta_{\tau(\ell_2)}) \in
[m^2,\omega)$ (so it is not equal to $\omega$). \newline
So we finish \footnote{see alternatively \scite{2.2}(1) + \scite{4.1}} 
proving $(*)_2$, hence $\otimes$.
\enddemo
\bigskip

\noindent
\underbar{Stage B}:  Like Stage A of \cite[III,4.4,p.164]{Sh:g}'s proof.
(So for $\alpha < \beta < \omega_2$, $\alpha$ does not appear in 
$\rho(\beta,\alpha)$).
\bigskip

\noindent
\underbar{Stage C}:  Defining the colouring:
\smallskip

Remember that ${\Cal S}^\alpha_\beta = \{ \delta < \aleph_\alpha:
\text{cf}(\delta) = \aleph_\beta\}$. \newline

For $\ell = 1,2$ choose $h_\ell:\omega_2 \rightarrow \omega_\ell$ such that
$S^\ell_\alpha = {\Cal S}^2_1 \cap h^{-1}_\ell(\{ \alpha\})$ is stationary for
each $\alpha < \omega_\ell$.  For $\alpha < \omega_2$, let 
$A_\alpha \subseteq \omega_1$ be such that no one is included in the 
union of finitely many others. \newline
\medskip

For $\alpha < \beta < \omega_2$, let $\ell = \ell_{\beta,\alpha}$ be
minimal such that

$$
d_{sq} \left( \rho_{h_1}(\beta,\alpha) \right) \in A_{\rho(\beta,\alpha)
(\ell)}
$$

\noindent
and lastly let

$$
c(\beta,\alpha) = c(\alpha,\beta) =: h_2 \biggl( (\rho(\beta,\alpha))
(\ell_{\beta,\alpha}) \biggr).
$$
\bigskip

\noindent
\underbar{Stage D}:  Proving that the colouring works:

So assume $n < \omega,\langle u_\alpha:\alpha < \omega_2 \rangle$ is a
sequence of pairwise disjoint subsets of $\omega_2$ of size $n$ and
$\gamma(*) < \omega_2$ and we should find $\alpha < \beta$ such that
$c \restriction (u_\alpha \times u_\beta)$ is constantly $\gamma(*)$. Without 
loss of generality
$\alpha < \beta \Rightarrow \text{ max}(u_\alpha) < \text{ min}(u_\beta)$
and $\text{min}(u_\alpha) > \alpha$
and let $E = \{\delta:\delta$ a limit ordinal $< \omega_2$ and
$(\forall \alpha)(\alpha < \delta \Rightarrow u_\alpha \subseteq \delta)\}$.
Clearly $E$ is a club of $\omega_2$.  For each $\delta \in E \cap 
{\Cal S}^2_1$, there is $\alpha^*_\delta < \delta$ such that

$$
\alpha \in [\alpha^*_\delta,\delta) \and \beta \in u_\delta \Rightarrow
\rho(\beta,\delta) \char 94 \langle \delta \rangle \, \trianglelefteq \,
\rho(\beta,\alpha).
$$
\medskip

\noindent
Also for $\delta \in {\Cal S}^2_1$ let

$$
\align
\varepsilon_\delta =: \text{ Min}\biggl\{ \varepsilon < \omega_1:&\zeta \in 
A_\delta \text{ but if } \alpha \in \dsize \bigcup_{\beta \in u_\delta}
\text{ Rang}(\rho(\beta,\delta)) \\
  &\text{(so } \alpha > \delta) \text{ then } \varepsilon \notin 
A_\alpha \biggr\}.
\endalign
$$
\medskip

\noindent
Note that $\varepsilon_\delta < \omega_1$ is well defined by the choice 
of $A_\alpha$'s.  So, by Fodor's lemma, for some $\zeta^* < \omega_1$
and $\alpha^* < \omega_2$ we have that

$$
W =: \{ \delta \in S^2_{\gamma(*)}:\alpha^*_\delta = \alpha^* \text{ and }
\varepsilon_\delta = \varepsilon^*\}
$$
\medskip

\noindent
is stationary.  Let $h$ be a strictly increasing function from $\omega_2$
into $W$ such that $\alpha^* < h(\delta)$.  By the demand on $\alpha^*$ (and
$W$)

$$
\alpha^* < \alpha < \delta \in W \and \beta \in u_\delta \Rightarrow
\rho(\beta,\delta) \char 94 \langle \delta \rangle \trianglelefteq 
\rho(\beta,\alpha). \tag"{$\bigoplus_0$}"
$$
\medskip

\noindent
Hence

$$
\align
\alpha^* < \alpha < \delta \in {\Cal S}^2_1 \and \beta \in 
u_{h(\delta)} \Rightarrow &\text{ Min}\{ \ell:\varepsilon^* \in 
A_{\rho(\beta,\alpha)(\ell)}\} =   \tag"{$\bigoplus_1$}"  \\
  &\text{ Min}\{ \ell:\rho(\beta,\delta)(\ell) = h(\delta)\}, 
\endalign
$$
\medskip

\noindent
hence

$$
\align
\alpha^* < \alpha < \delta \in {\Cal S}^2_1 \and &\beta \in 
u_{h(\delta)} \Rightarrow    \tag"{$\bigoplus_2$}"   \\
  &h_2 \biggl( \rho(\beta,\delta)\biggl[ \text{Min}\{ \ell:\varepsilon^* \in
A_{\rho(\beta,\delta)(\ell)}\}\biggr] \biggr) = \gamma(*). 
\endalign
$$
\medskip

\noindent
Let

$$
\align
E_0 =: \biggl\{ \delta < \omega_2:&\delta \text{ a limit ordinal, } 
\delta \in E \text{ and}\\
  &\alpha < \delta \Rightarrow h(\alpha) < \delta \text{ (hence }
\sup(u_{h(\alpha)}) < \delta) \biggr\}.
\endalign
$$
\medskip

\noindent
For each $\delta \in {\Cal S}^2_1$ there is $\alpha^{**}_\delta < \delta$
such that $\alpha^{**}_\delta > \alpha^*$ and

$$
\alpha \in [\alpha^{**}_\delta,\delta) \and \beta \in u_{h(\delta)}
\Rightarrow \rho(\beta,\delta) \char 94 \langle \delta \rangle
\trianglelefteq \rho(\beta,\alpha).
$$
\medskip

\noindent
For each $\gamma < \omega_1,\delta \mapsto \alpha^{**}_\delta$ is a regressive
function on $S^1_\gamma$, hence for some \newline
$\alpha^{**}(\gamma) < \delta$ the set
$S'_\gamma =: \{ \delta \in S^1_\gamma \cap E_0:\alpha^{**}_\delta =
\alpha^{**}(\gamma)\}$ is stationary.

Let $\alpha^{**} = \text{ sup}\{ \alpha^{**}
(\gamma)+1:\gamma < \omega_1\}$ and note that $\alpha^{**} < \omega_2$.  Let

$$
E_1 =: \{ \delta < \omega_2:\text{for every } \gamma < \omega_1,\delta =
\text{ sup}(S'_\gamma \cap \delta) \text{ and } \delta > \alpha^{**}\},
$$
\medskip

\noindent
and note that $E_1$ is a club of $\aleph_2$ (and as $S'_\gamma \subseteq E_0$ 
clearly $E_1 \subseteq E_0$) and choose $\delta^* \in E_1 \cap 
S^2_{\gamma(*)}$.  Then by induction on $i < \omega_1$ choose an ordinal 
$\zeta_i$ such that $\langle \zeta_i:i < \omega_1 \rangle$ is strictly 
increasing with limit $\delta^*$ and 
$\zeta_i \in S'_i \backslash (\alpha^{**} + 1)$.
We know that $\alpha < \zeta_i \Rightarrow u_\alpha \subseteq
\zeta_i$ and $\alpha < \text{ min}(u_\alpha)$, hence for every
$\alpha_i < \zeta_i$ large enough $(\forall \beta \in u_{\alpha_i})
(\rho(\delta^*,\zeta_i) \char 94 (\zeta_i) \trianglelefteq \rho(\delta^*,
\beta))$.

Choose such $\alpha_i \in (\dsize \bigcup_{j < i} \zeta_j,\zeta_i)$.
Lastly for $i < \omega_1$ choose $\beta_i \in E \cap S'_i$ with $\beta_i > 
\delta^*$.  Now for each $i < \omega_1$ for some $\xi(i) < \delta^*$,

$$
\alpha \in (\xi(i),\delta^*) \and \beta \in u_{h(\beta_i)} \Rightarrow
\rho(\beta,\delta^*) \char 94 \langle \delta^* \rangle \, \trianglelefteq \,
\rho(\beta,\alpha). \tag"{$\bigoplus_3$}"
$$

\medskip
\noindent
As $\delta^* = \dsize \bigcup_{i < \omega_1} \zeta_i$, without loss of
generality $\xi(i) = \zeta_{j(i)}$, and $j(i)$ is (strictly) increasing
with $i$ and let
$A =: \{ \varepsilon < \omega_1:\varepsilon \text{ a limit ordinal and }
(\forall i < \varepsilon)(j(i) < \varepsilon)\}$.  Clearly $A$ is a club of
$\omega_1$.  Now putting all of this together we have:
\medskip
\roster
\item "{$(*)_1$}"  if $i(0) < i(1)$ are in $A,\alpha \in u_{\alpha_{i(1)}},
\beta \in u_{h(\beta_{i(0)})}$ then \newline
$\rho(\beta,\alpha) = \rho(\beta,\delta^*) \char 94 \rho(\delta^*,\alpha)$.
\newline
[Why?  As $j(i(0)) < i(1)$, see $\bigoplus_3$].
\item "{$(*)_2$}"  if $i < \omega_1$ then $\beta \in u_{h(\beta_i)} 
\Rightarrow
i \in \text{ Rang}(\rho_{h_1}(\beta,\delta^*))$ (witnessed by $\beta_i$ which
belongs to this set by $\bigoplus_1$).
\item "{$(*)_3$}"  if $i < \omega_1$ then $\alpha \in u_{\alpha_i} \Rightarrow
i \in \text{ Rang}(\rho_{h_1}(\delta^*,\alpha))$ (witnessed by $\zeta_i$
which belongs to this set by the choice of $\alpha_i$)
\item "{$(*)_4$}"  if $i < \omega_1$ and $\beta \in u_{h(\beta_i)}$ then
$\ell = \text{ Min}\{\ell:\zeta^* \in A_{\rho(\beta,\delta^*)(\ell)}\}$ is
well defined and $h_2(\rho(\beta,\delta^*)(\ell)) = \gamma(*)$. \newline
[Why?  By $\bigoplus_2$].
\endroster
\medskip

\noindent
Now let $\nu_i$, for $i < \omega_1$, be the concatanation of $\{ \rho
(\beta,\delta^*):\beta \in u_{\beta_i}\}$ and $\rho_i$ be the concatanation
of $\{\rho(\delta^*,\alpha):\alpha \in u_{\alpha_i}\}$.  So we can apply
$\otimes$ of Stage A to \newline
$\langle \nu_i,\rho_i:i < \omega_1 \rangle$ and
$\gamma^*$ (its assumptions hold by $(*)_1 + (*)_2 + (*)_3$) and get that
for some $i < j < \omega_1$ we have $d_{\text{sq}}(\nu' \char 94 \rho') 
= \zeta^*$
whenever $\nu'$ is a subsequence of $\nu_i,\rho'$ a subsequence of $\rho_j$
such that $i \in \text{ Rang}(\nu'),j \in \text{ Rang}(\rho')$.  Now for
$\beta \in u_{h(\beta_i)}$, \newline
$\alpha \in u_{\alpha_j}$ we have
\medskip
\roster
\item "{{}}"  $\rho(\beta,\alpha) = \rho(\beta,\delta^*) \char 94 
\rho(\delta^*,\alpha) \text{ (see } (*)_1) \text{ and}$
\medskip
\item "{{}}"  $\rho(\beta,\delta^*)$ is O.K. as $\nu'$.
{\roster
\itemitem { {} }
[Why?  Because it is a subsequence of $\nu_i$ (see the choice of $\nu_i$) and
$i$ belongs to Rang$(\rho(\beta,\delta^*))$ by $(*)_2$] and
\endroster}
\item "{{}}"  $\rho(\delta^*,\alpha)$ is O.K. as $\rho'$
{\roster
\itemitem { {} }
[Why?  Because $\rho(\delta^*,\alpha)$ is a subsequence of $\rho_j$ by the
choice of $\rho_j$ and $j$ belongs to Rang$(\rho(\delta^*,\alpha))$ by
$(*)_3$].
\endroster}
\endroster
\medskip

\noindent
Now by $(*)_4$ the colour $c(\beta,\alpha)$ is $\gamma(*)$ as 
required and get the desired conclusion. \newline
${}$ \hfill$\square_{\scite{1.1}}$
\bigskip

\remark{Remark}  Can we get $Pr_1(\lambda^{+2},\lambda^{+2},\lambda^{+2},
\lambda)$ for $\lambda$ regulars by the above proof? If $\lambda = 
\lambda^{< \lambda}$ the
same proof works (now $\text{Dom}(d_{sq}) = {}^{\omega >}(\lambda^+)$ and
$\nu_\alpha,\rho_\alpha \in {}^{\lambda >}(\lambda^+))$. \newline
See more in \S2.
\endremark
\newpage

\head {\S2 Larger Cardinals: The implicit properties} \endhead
\resetall
\bigskip

\noindent
More generally (than in the remark at the end of \S1):
\definition{\stag{2.1} Definition}  1) $Pr_6(\lambda,\lambda,\theta,\sigma)$ 
means
that there is $d:{}^{\omega >} \lambda \rightarrow \theta$ such that:
\newline
\underbar{if} $\langle (u_\alpha,v_\alpha):\alpha < \lambda \rangle$
satisfies

$$
u_\alpha \subseteq {}^{\omega >}\lambda,v_\alpha \subseteq {}^{\omega >}
\lambda,
$$

$$
|u_\alpha \cup v_\alpha| < \sigma,
$$

$$
\nu \in u_\alpha \cup v_\alpha \Rightarrow \alpha \in \text{ Rang}(\nu),
$$
\medskip

\noindent
and $\gamma < \theta$ and $E$ a club of $\lambda$ \underbar{then} for
some $\alpha < \beta$ from $E$ we have

$$
\nu \in u_\alpha \and \rho \in v_\beta \Rightarrow d(\nu \char 94 \rho)
= \gamma.
$$
\medskip

\noindent
2)  $Pr^6_S(\lambda,\lambda,\theta,\sigma)$ is defined similarly but
$\alpha < \beta$ are required to be in $E \cap S$.  $Pr^6_\tau(\lambda,
\lambda,\theta,\sigma)$ means ``for some stationary $S \subseteq
\{ \delta < \lambda:\text{cf}(\delta) \ge \tau\}$ we have
$Pr^6_S(\lambda,\lambda,\theta,\sigma)$".  If $\tau$ is 
omitted, we mean $\tau = \sigma$.  Lastly $Pr^6_{\text{nacc}}(\lambda,
\lambda,\theta,\sigma)$ is defined similarly but demanding $\alpha,\beta \in
\text{ nacc}(E)$ and $Pr^-_6(\lambda,\lambda,\theta,\sigma)$ is defined
similarly but $E = \lambda$.
\enddefinition  
\bigskip

\proclaim{\stag{2.2} Lemma}  0) If $Pr_6(\lambda,\lambda,\theta,\sigma)$ and
$\theta_1 \le \theta$ and $\sigma_1 \le \sigma$ \underbar{then} $Pr_6(\lambda,
\lambda,\theta_1,\sigma_1)$ (and similar monotonicity properties for 
Definition \scite{2.1}(2)).
Without loss of generality $u_\alpha = v_\alpha$ in Definition \scite{2.1}.
\newline
1) If $Pr_6(\lambda^+,\lambda^+,\lambda^+,\lambda)$, \underbar{then}
$Pr_1(\lambda^{+2},\lambda^{+2},\lambda^{+2},\lambda)$. \newline
2)  If $Pr_6(\lambda^+,\lambda^+,\theta,\sigma)$, so $\theta \le \lambda^+$ 
\underbar{then}
$Pr_1(\lambda^{+2},\lambda^{+2},\lambda^{+2},\sigma)$ provided that
\medskip
\roster
\item "{$(*)$}"  there is a sequence $\bar A = \langle A_\alpha:\alpha < 
\lambda^{++} \rangle$ of subsets of $\theta$ such that for every 
$\alpha \in u \subseteq \lambda^{++}$ with $u$ of cardinality $< \sigma$,
we have

$$
A_\alpha \backslash \cup \{A_\beta:\beta \in u,\beta \ne \alpha\} \ne 
\emptyset.
$$
\endroster
\medskip

\noindent
3) If $\lambda$ is regular and $\lambda = \lambda^{< \lambda}$ \underbar{then}
$Pr_6(\lambda^+,\lambda^+,\lambda^+,\lambda)$. \newline
4)  In \cite[III,4.7]{Sh:g} we can change the assumption accordingly.
\endproclaim
\bigskip

\demo{Proof}  0) Clear. \newline
1) By part (2) choosing $\theta = \lambda^+,\sigma = \lambda$ as $(*)$ 
holds as $\lambda^+$ is regular (so e.g. choose by induction on $\alpha <
\lambda^{++},A_\alpha \subseteq \lambda^+$ see unbounded non-stationary with
$\beta < \alpha \Rightarrow |A_\alpha \cap A_\alpha| \le \lambda$.\newline
2)  Like the proof for $\aleph_2$, only now $\{ \delta < \lambda^{++}:
\text{cf}
(\delta) = \lambda^+\}$ plays the role of ${\Cal S}^2_1$ and let
$h_1:\lambda^{++} \rightarrow \lambda^+$ and $h_2:\lambda^{++}
\rightarrow \lambda^{++}$ be such that for every $\gamma < \lambda^{+ \ell}$
and $\ell \in \{ 1,2\}$ the set $S^\ell_\gamma = \{ \alpha < \lambda^{+2}:
\text{cf}(\alpha) = \lambda^+ \text{ and } h_\ell(\alpha) = \gamma\}$ is
stationary.  Finally, if $dq$
exemplifies $Pr_6(\lambda^+,\lambda^+,\theta,\sigma)$, then in defining
$c$ for a given 
$\alpha < \beta$, let $\ell_{\alpha,\beta}$ be the minimal $\ell$ such that
$dq(\rho_{h_1}(\alpha,\beta))$ belongs to 
$A_{\rho_{h_1}(\alpha,\beta)(\ell)}$ and 
let $c(\beta,\alpha) = c(\alpha,\beta) = h_2\left( \rho(\beta,\alpha)
(\ell_{\beta,\alpha}) \right)$.  Then in stage D without loss of generality
$|u_\alpha| = \sigma < \lambda$ for $\alpha < \lambda^+$ and continue as 
there, but after the definition of $E_1$ we do not choose $\zeta_i,\alpha_i$
instead we first continue choosing $\beta_i,\xi_i$ for $i < \lambda^+$ as
there as without loss of generality $\delta^* = \dsize \bigcup_{i < \lambda^+}
\xi(i)$.  Only then we choose by induction on $i < \lambda^+$ the ordinal
$\zeta_i$ such that: $\zeta_i \in S'_i \backslash (\alpha^{**} +1),\zeta_i
> \sup[\{\xi(j):j \le i\} \cup \{\zeta_j:j < i\}]$ and then choose
$\alpha_i < \zeta_i$ large enough (so no need of the club $A$ of 
$\lambda^+$). \newline
3) As in the proof of \scite{1.1}, Stage A. \newline
4) Combine the proofs here and there (and not used). 
\hfill$\square_{\scite{2.2}}$
\enddemo
\bigskip

\noindent
This leaves some problems on $Pr_1$ open; e.g. \newline
\subhead {\stag{2.3} Question} \endsubhead  1) If $\lambda > \aleph_0$ is 
inaccessible,
do we have $Pr_1(\lambda^+,\lambda^+,\lambda^+,\lambda)$
(rather than with $\sigma < \lambda$)? \newline
2) If $\mu > \aleph_0$ is regular (singular) and $\lambda = \mu^+$, do we have
$Pr_1(\lambda^+,\lambda^+,\lambda^+,\mu)$?
\newline
[clearly, yes, for the weaker version: $c$ a symmetric two place function
from $\lambda^+$ to $\lambda^+$
such that for every $\gamma < \lambda^+$ and pairwise disjoint
$\langle u_\alpha:\alpha < \lambda^+ \rangle$ with 
$u_\alpha \in [\lambda^+]^{< \lambda}$ we have

$$
(\exists \alpha < \beta)\forall i \in u_\alpha \, \forall j \in u_\beta
\biggl( \gamma \in \text{ Rang }\rho_c(j,i) \biggr)].
$$
\medskip

\noindent
See more in \S4.  Remember that we know $Pr_1(\lambda^{+2},\lambda^{+2},
\lambda^{+2},\sigma)$ for $\sigma < \lambda$.
\bigskip

\proclaim{\stag{2.4} Claim}  Assume $\mu$ is singular, $\lambda = \mu^+$,
$2^\kappa > \mu > \kappa = \kappa^\theta,\theta = \text{cf}(\theta) 
\ge \sigma$ and $Pr_6(\theta,\theta,\theta,\sigma)$.
\underbar{Then} $Pr_1(\mu^+,\mu^+,\theta,\sigma)$. 
\endproclaim
\bigskip

\demo{Proof}  Let $\bar e = \langle e_\alpha:\alpha < \lambda \rangle$ be
a club system, $S \subseteq \{ \delta < \mu^+:cf(\delta) = \theta\}$
stationary such that $\lambda \notin \text{ id}^a(\bar e \restriction S)$
and $\alpha \in e_\delta \Rightarrow \text{ cf}(\alpha) \ne \theta$ and

$$
\delta = \sup(\delta \cap S) \and \chi < \mu 
\Rightarrow \delta = \text{ sup} \left(
\{ \alpha \in e_\delta:cf(\alpha) > \chi + \sigma^+,\text{ so }
\alpha \in \text{ nacc}(e_\delta)\} \right)
$$

\noindent
and $\alpha \in e_\beta \cap S \Rightarrow e_\alpha \subseteq e_\beta$
(exists by \cite[2.10]{Sh:365}). Let \newline
$\bar f = \langle f_\alpha:\alpha < \theta \rangle,f_\alpha:\mu^+
\rightarrow \kappa$ such that every partial function $g$ from $\mu^+$ to
$\kappa$ (really $\sigma$ suffice) of cardinality $\le \theta$ is included
in some $f_\alpha$ (exist by \cite{EK} or see \newline
\cite[AP1.7]{Sh:g}).
\medskip

So for some $f = f_{\alpha(*)}$ we have
\medskip
\roster
\item "{$(*)$}"  for every club $E$ of $\mu^+$ for some $\delta \in S$
we have:
{\roster
\itemitem{ (a) }  $e_\delta \subseteq E$
\itemitem{ (b) }  if $\chi < \mu$ and $\gamma < \theta$ then \newline
$\delta = \text{ sup}(\{ \alpha \in \text{ nacc}(e_\delta):
f(\alpha) = \gamma \text{ and } \text{cf}(\alpha) > \chi\})$.
\endroster}
\endroster
\medskip

\noindent
This actually proves $\text{id}_p(\bar e \restriction S)$ is not weakly
$\theta^+$-saturated. \newline
The rest is by combining the trick of \cite[III,\S4]{Sh:g} (using first
$\delta(*) \in S$ then some suitable $\alpha \in \text{ nacc}(e_{\delta(*)})$)
and the proof for $\aleph_2$. \hfill$\square_{\scite{2.4}}$
\enddemo
\bigskip

\subhead {\stag{2.5} Fact} \endsubhead  $Pr_1(\lambda^+,\lambda^+,\theta,
\text{cf}(\lambda))$ implies $Pr^6(\lambda^+,\lambda^+,\theta,
\text{cf}(\lambda))$.
\bigskip

\remark{Remark}  This is not totally immediate as in $Pr_1$ the sets are
required to be pairwise disjoint.
\endremark
\bigskip

\demo{Proof}  Let $\kappa = \text{ cf}(\lambda)$ and $f_\alpha \in 
{}^\kappa \lambda$ for $\alpha < \lambda^+$
be such that $\alpha < \beta \Rightarrow f_\alpha <^*_{J^{bd}_\kappa}f_\beta$.
Let $d:[\lambda^+]^2 \rightarrow \theta$ exemplifies $Pr_1(\lambda^+,
\lambda^+,\theta,\text{cf}(\lambda))$.  Let $c:\kappa \rightarrow \kappa$ be
such that for every $\gamma < \kappa$ for undoubtedly many $\beta < \kappa$
we have $c(\beta) = \gamma$.  For $\nu \in {}^{\omega >}(\lambda^+)$ 
we define $d^*_{sq}(\nu)$ as follows.
\medskip

\noindent
If $\ell g(\nu) \le 1$ or $\nu$ is constant, then $d^*_{\text{sq}}(\nu) = 0$.
So assume $\ell g(\nu) \ge 2$ and $\nu$ is not constant.

$$
\text{For }\alpha < \beta < \lambda^+ \text{ let }
{\bold s}(\beta,\alpha) = {\bold s}(\alpha,\beta) = \text{ sup}\{i+1:i 
< \kappa \text{ and } f_\alpha(i) \ge f_\beta(i)\},
$$

$$
{\bold s}(\alpha,\alpha) = 0,
$$

$$
{\bold s}(\nu) = \text{ max}\{ {\bold s}(\nu(\ell),\nu(k)):\ell,
k < \ell g(\nu) \text{ (so } {\bold s} \text{ is symmetric)}\},
$$

$$
a(\nu) = \{(\ell,k):{\bold s}(\nu(\ell),\nu(k)) = \bold s(\nu) \text{ and }
\ell < k < \ell g(\nu)\}.
$$
\medskip

\noindent
As $\ell g(\nu) \ge 2$ and $\nu$ is not constant, clearly $a(\nu) \ne 
\emptyset$ and $a(\nu)$ is finite, so let
$(\ell_\nu,k_\nu)$ be the first pair from $a(\nu)$ in lexicographical
ordering. \newline
\smallskip

\noindent
Lastly $d^*_{\text{sq}}(\nu) = c 
\biggl( d(\{\nu(\ell_\nu),\nu(k_\nu)\}) \biggr)$.
\bigskip

Now we are given $\gamma < \theta$, stationary 
$S \subseteq \{ \delta < \lambda^+:\text{cf}(\delta) \ge \text{ cf}
(\lambda)\}$,
$\langle u_\alpha:\alpha < \lambda^+ \rangle$ (remember \scite{2.2}(0)),
$|u_\alpha| < \text{ cf}
(\lambda),u_\alpha \subseteq {}^{\omega >} \lambda$ such that 
$\alpha \in \cap
\{ \text{Rang}(\nu):\nu \in u_\alpha\}$. \newline
Let $u'_\alpha = \cup \{\text{Rang}(\nu):\nu \in u_\alpha\}$ and 
without loss of generality for some stationary $S' \subseteq S$ and 
$\gamma_0,\beta^*$ we have $\alpha \in S' \Rightarrow \gamma_0 =
\text{min}\{\gamma + 1:\text{if } \beta_1 < \beta_2$
are in $u'_\alpha$ then $f_{\beta_1} \restriction [\gamma,\text{cf}(\lambda))
< f_{\beta_2} \restriction [\gamma,\text{cf}(\lambda))\} < \kappa$ and 
$\sup (\cup \{ u'_\alpha \cap \alpha:\alpha \in S'\}) < \beta^* < \lambda^+$.
Now for some $\gamma_1 \in (\gamma_0,\text{cf}(\lambda))$ and stationary 
$S_0,S_1 \subseteq S'$ and $\gamma^* < \lambda$ we have

$$
\beta \in u'_\alpha \and \alpha \in S_0 \Rightarrow f_\beta(\gamma_1) <
\gamma^*,
$$

$$
\beta \in u'_\alpha \and \alpha \in S_1 \Rightarrow f_\beta(\gamma_1) >
\gamma^*.
$$
\medskip

\noindent
Let $\{ \alpha^\ell_\zeta:\zeta < \lambda\}$ enumerate some unbounded
$S'_\ell \subseteq S_\ell$ in increasing order such that
$\zeta < \xi \Rightarrow \sup(u_{\alpha^0_\zeta} \cup u_{\alpha^1_\zeta})
< \text{ min}(u_{\alpha^0_\xi} \cup u_{\alpha^1_\xi})$.

Lastly apply the choice of $d$. \hfill$\square_{\scite{2.5}}$
\enddemo
\newpage

\head {\S3 Guessing Clubs Revisited} \endhead
\resetall
\bigskip

\proclaim{\stag{3.1} Claim}  Assume $\lambda = \mu^+$, and \newline
$S \subseteq \{ \delta < \lambda^+:\text{cf}(\delta) = \lambda$ and 
$\delta$ is
divisible by $\lambda^2\}$ is stationary. \newline
1)  There is a strict club system 
$\bar C = \langle C_\delta:\delta \in S \rangle$ such that 
$\lambda^+ \notin \text{ id}^p(\bar C)$ and \newline
$[\alpha \in \text{ nacc}
(C_\delta) \Rightarrow cf(\alpha) = \lambda]$; moreover, there are
$h_\delta:C_\delta \rightarrow \mu$ such that for every club $E$ of
$\lambda^+$, for stationarily many $\delta \in S$,

$$
\dsize \bigwedge_{\zeta < \mu} \delta = \sup [h^{-1}_\delta(\{ \zeta \})
\cap E \cap \text{ nacc}(C_\delta)].
$$
\medskip

\noindent
2) If $\bar C$ is a strict $S$-system, $\lambda^+ \notin \text{ id}^p
(\bar C,\bar J),
J_\delta$ a $\lambda$-complete ideal on $C_\delta$ extending
$J^{bd}_{C_\delta} + \text{ acc}(C_\delta)$ (with $S,\mu$ as above) 
\underbar{then} the parallel result holds for some \newline
$\bar h = \langle h_\delta:\delta \in S \rangle$ where $h_\delta$ is a
function from $C_\delta$ to $\mu$, i.e. we have for every club $E$ of
$\lambda^+$, for stationarily many $\delta \in S \cap \text{ acc}(E)$ for
every $\gamma < \mu$ the set \newline
$\{ \alpha \in C_\delta:h_\delta(\alpha) =
\gamma \text{ and } \alpha \in E\}$ is $\ne \emptyset \text{ mod } J_\delta$.
\endproclaim
\bigskip

\remark{\stag{3.1A} Remark}  1) This improves \cite[3.1]{Sh:413}. \newline
2) Of course, we can strengthen (1) to:

$$
\{ \alpha \in C_\delta:h_\delta(\alpha) = \gamma \text{ and } \alpha \in E
\text{ and } \alpha \in \text{ nacc}(C_\delta) \text{ and } \sup(\alpha \cap
C_\delta) \in E\}.
$$
\medskip

\noindent
E.g. for every thin enough club $E$ of $\lambda,\bar C^E$ will serve where:
$C^E_\delta = C_\delta \cap E$ if $\delta \in \text{ acc}(E)$ and $C^E_\delta
= C_\delta$, otherwise. \newline
For \scite{3.1}(2) we get slightly less: for some club $E^*:\{\alpha \in
C_\delta:h_\delta(\alpha) = \gamma \text{ and } \alpha \in E \text{ and }
\alpha \in \text{ nacc}(C_\delta) \text{ and } \sup(\alpha \cap C_\delta
\cap E^*) \in E\}$.
\endremark
\bigskip

\demo{Proof}  1) Let $\langle C_\delta:\delta \in S \rangle$ be such that
$\lambda^+ \notin \text{ id}^p(\bar C)$ and \newline
$[\alpha \in \text{ nacc}(C_\delta)
\Rightarrow cf(\delta) = \lambda]$ (such a sequence exists by 
\cite[2.4(3)]{Sh:365}).
Let $J_\delta = J^{bd}_{C_\delta} + \text{ acc}(C_\delta)$.  Now apply
part (2).  \newline
2) For each $\delta \in S$ let $\langle A^\alpha_\delta:\alpha \in C_\delta
\rangle$ be a sequence of distinct non-empty subsets of $\mu$ to be 
chosen later.  By induction on $\zeta < \lambda$ we try to define
$E_\zeta,\langle Y^\zeta_\alpha:\alpha \in S \rangle$, \newline
$\langle Z^\zeta_{\alpha,\gamma}:\alpha \in \zeta
\text{ and } \gamma < \mu \rangle$ such that

$$
E_\zeta \text{ is a club of } \lambda^+, \text{ decreasing in } \zeta,
$$
\medskip

\noindent
for $\gamma < \mu$,

$$
Z^\zeta_{\delta,\gamma} = \{ \alpha:\alpha \in E_\zeta \cap
\text{ nacc}(C_\delta) \text{ and } \gamma \in A^\alpha_\delta\},
$$

$$
Y^\zeta_\delta = \{ \gamma < \mu:Z^\zeta_{\delta,\gamma} \ne \emptyset
\text{ mod } J_\delta\}.
$$
\medskip

\noindent
$E_{\zeta + 1}$ is such that

$$
\align
\biggl\{ \delta \in S:&Y^\zeta_\delta = Y^{\zeta + 1}_\delta 
\text{ and } \delta \in \text{ nacc}(E_{\zeta + 1}) \\
  &\text{and } E_{\zeta +1} \cap \text{ nacc}(C_\delta) \notin J_\delta 
\biggr\} \text{ is not stationary}.
\endalign
$$
\medskip

\noindent
If we succeed to define $E_\zeta$, for each $\zeta < \lambda$, then
$E =: \dsize \bigcap_{\zeta < \lambda} E_\zeta$ is a club of $\lambda^+$, 
and since
$\lambda^+ \notin \text{ id}^p(\bar C)$, we can choose 
$\delta \in S$ such that
$\delta = \sup(E \cap \text{ nacc }C_\delta)$ and $E \cap \text{ nacc}
(C_\delta) \ne \emptyset \text{ mod } J_\delta$.  Then as 
$\dsize \bigcup_{\gamma < \mu} Z^\zeta_{\delta,\gamma} \supseteq E \cap 
\text{ nacc}(C_\delta)$ for each $\zeta < \lambda$ necessarily (by the
requirement on $J_\delta$) for some $\gamma < \mu,Z^\zeta_{\delta,\gamma}
\ne \emptyset \text{ mod }J_\delta$, hence
$Y^\zeta_\delta \ne \emptyset$ so that 
$\langle Y^\zeta_\delta:\zeta < \lambda \rangle$ is a strictly 
decreasing sequence of subsets of
$\mu$, and since
$\mu < \text{ cf}(\mu^+) = \text{ cf}(\lambda)$, we have a contradiction.
So necessarily we will be stuck (say) for $\zeta(*) < \lambda$. \newline
We still have the freedom of choosing $A^\alpha_\delta$ for
$\alpha \in C_\delta$.
\enddemo
\bigskip

\noindent
\underbar{Case 1}:  $\mu$ regular. \newline
By induction on $\alpha \in C_\delta$ we can choose sets $A^\alpha_\delta$
such that
\medskip
\roster
%
\item "{(i)}"  $A^\alpha_\delta \subseteq \mu,|A^\alpha_\delta| = \mu,
\langle A^\alpha_\delta:\alpha \in C_\delta,\text{otp}(\alpha \cap C_\delta) <
\mu \rangle$ are pairwise disjoint,
\item "{(ii)}"  for $\beta \in C_\delta \cap \alpha$,
$A^\alpha_\delta \cap A^\beta_\delta$ is bounded in $\mu$,
\item "{(iii)}" if $\mu > \aleph_0$ then $A^\alpha_\delta$ is non-stationary
(just to clarify their choice).
\endroster
\medskip

\noindent
There is no problem to carry the induction. \newline
We shall prove later that
\medskip
\roster
\item "{$(*)$}"  if $E$ is a club of $\lambda^+,\delta \in S \cap 
\text{ acc}(E)$ and
$\delta = \sup(E \cap \text{ nacc }C_\delta)$ and \newline
$E \cap \text{ nacc}(C_\delta) \ne \emptyset \text{ mod } J_\delta$
\underbar{then}
\endroster
\medskip

\roster
\item "{$(**)_\delta$}"  for some $\alpha_\delta \in E \cap
\text{ nacc}(C_\delta)$, the following set $B_\delta$ is unbounded in $\mu$,
where

$$
\align
B_\delta = \biggl\{ \gamma < \mu:&\{\beta:\beta \in
E \cap \text{ nacc}(C_\delta) \text{ and } \beta \ne 
\alpha_\delta \\
  &\text{ and } \gamma = \sup(A^{\alpha_\delta}_\delta \cap A^\beta_\delta)\}
\ne \emptyset \text{ mod }J_\delta \biggr\}.
\endalign
$$
\endroster
\medskip

\noindent
Choose the minimal such that $\alpha_\delta = \alpha^E_\delta$ 
(for other $\delta$'s it does not matter, i.e. for those for which
$\delta > \sup(E \cap \text{ nacc}(C_\delta))$ or
$E_{\zeta(*)} \cap \text{ nacc}(C_\delta) \in J_\delta$). \newline
Clearly if $E' \supseteq E''$ and $\alpha^{E'}_\delta,\alpha^{E''}
_\delta$ are defined then $\alpha^{E'}_\delta \le \alpha^{E''}_\delta$. 
\newline
Now for any club $E^* \subseteq E_{\zeta(*)}$ of $\lambda^+$, for
$\delta \in S \cap \text{ acc}(E_{\zeta(*)})$ we define \newline
$h^{E^*}_\delta:C_\delta \rightarrow \mu$ by letting
$h^{E^*}_\delta(\beta) = \text{ otp}(B_\delta \cap \sup
(A^{\alpha_\delta}_\delta \cap A^\beta_\delta))$ for $\beta \in C_\delta
\backslash \{ \alpha_\delta\}$ and 
$h^{E^*}_\delta(\alpha_\delta) = 0$. \newline
\medskip

\noindent
Now for any club $E$ of $\lambda^+$ for stationarily many $\delta \in S
\cap \text{ acc}(E^* \cap E)$, we have

$$
\biggl\{
\gamma < \mu:\{ \alpha:\alpha \in E^* \cap E \cap E_{\zeta(*)} \cap
\text{ nacc}(C_\delta) \text{ and } \gamma \in A^\alpha_\delta\} \ne
\emptyset \text{ mod } J_\delta \biggr\} = Y^{\zeta(*)}_\delta
$$
\medskip

\noindent
(this holds by the choice of $\zeta(*)$).  Let the set of such $\delta \in
S \cap \text{ acc}(E^* \cap E)$ be called $Z^{E^*}_E$.  Now for each
$\delta \in Z^{E^*}_E$, the set

$$
\align
B_\delta[E,E^*] =: \biggl\{ \gamma < \mu:&\{ \beta:\beta \in E
\cap E^* \cap E_{\zeta(*)} \cap \text{ nacc}(C_\delta) \\
  &\text{ and } \beta \ne \alpha^{E^*}_\delta \text{ and }
  \gamma = \sup(A^{\alpha_\delta}_\delta \cap A^\beta_\delta)\} \ne \emptyset
\text{ mod }J_\delta \biggr\}
\endalign
$$
\medskip

\noindent
is necessarily
unbounded in $\mu$.  So in the same way we have gotten $E_{\zeta(*)}$ we
can find club $E^* \subseteq E_{\zeta(*)}$ such that for any club $E$ of
$\lambda^+$, for stationarily many $\delta \in Z^{E^*}_E$ we have
$B_\delta[E,E_{\zeta(*)}] = B_\delta[E^*,E_{\zeta(*)}]$ and
$\alpha^E_\delta = \alpha^{E^*}_\delta$ (note the minimality in the choice
of $\alpha^E_\delta$ so it can change $\le \lambda + 1$ times; more
elaborately if $\langle E^*_\zeta:\zeta < \lambda \rangle$ is a decreasing 
sequence of clubs and $\delta \in
Z^{E^*}_{E^*}$, where $E^* = \dsize \bigcap_{\zeta < \lambda} E^*_\zeta$, then
$\langle \alpha^{E^*_\zeta}_\delta:\zeta < \lambda \rangle$ is increasing and
bounded in $C_\delta$ (by $\alpha^{E^*}_\delta$), hence is 
eventually constant).  
Define $h_\delta:C_\delta \rightarrow \mu$ by 
$h_\delta(\beta) = \text{ otp} \left (B_\delta[E^*,E_{\zeta(*)}]
\cap \sup(A^{\alpha_\delta}_\delta \cap A^\beta_\delta) \right)$
if $\beta \ne \alpha_\delta$ and $h_\delta(\beta) = 0$ if 
$\beta = \alpha_\delta$. \newline
\medskip

\noindent
Why does $(*)$ hold?  \newline
If not, let $B = E_{\zeta(*)} \cap \text{ nacc}(C_\delta)$, so
$\text{otp}(B) = \lambda = \mu^+$ and $B \ne \emptyset \text{ mod } 
J_\delta$, 
so for every $\alpha \in B$ we can find
$\varepsilon_\alpha < \mu$ and $Y_{\alpha,\varepsilon} \in J_\delta$ (for
$\varepsilon < \mu)$ such that if $\xi \in B \backslash Y_{\alpha,\varepsilon}
\backslash \{ \alpha \}$ and $\varepsilon \in [\varepsilon_\alpha,\mu)$ then
$\sup(A_\delta^\alpha \cap A_\delta^\xi) \ne \varepsilon$. 
Now let $Y_\alpha =: \cup\{ Y_{\alpha,\varepsilon}:\varepsilon \in
[\varepsilon_\alpha,\mu)\} \cup \{ \alpha + 1 \}$ and note that
$Y_\alpha \in J_\delta$.  So for some $\varepsilon^* < \mu$,
$B_1 =: \{ \alpha \in B:\varepsilon_\alpha =
\varepsilon^*\}$ is $\ne \emptyset \text{ mod }J_\delta$.
For each $\alpha \in B_1$ choose
$\gamma_\alpha \in A^\delta_\alpha \backslash (\varepsilon^* + 1)$
(remember $|A^\delta_\alpha| = \mu$).  So for some $\gamma^* < \mu$ the set
$B_2 =: \{ \alpha \in B_1:\gamma_\alpha = \gamma^*\}$ is $\ne \emptyset
\text{ mod }J_\delta$.  Let $\alpha^* = \text{ Min}(B_2)$, and for $\gamma
\in [\gamma^*,\mu)$ we define \newline
$B_{\zeta,\gamma} = \{ \alpha \in B_2:\gamma
= \sup(A_\delta^{\alpha^*} \cap A_\delta^\alpha\}$.  So clearly
$B_2 = \cup\{ B_{\zeta,\gamma}:\gamma^* \le \gamma < \mu\}$, hence for some
$\gamma^{**} \in [\gamma^*,\mu)$ we have $B_{\zeta,\gamma^{**}} \ne \emptyset
\text{ mod }J_\delta$, hence $\gamma^{**}$ contradicts the choice of
$\varepsilon_{\alpha^*} = \varepsilon^*$.
\bigskip

\noindent
\underbar{Case 2}:  $\mu$ singular. \newline
Let $\kappa = cf(\mu)$, so by \cite[II,\S1]{Sh:g} we can find an increasing
sequence $\langle \lambda_i:i < \kappa \rangle$ of regular cardinals
$> \kappa$ with limit $\mu$ such that $\lambda = \mu^+ = \text{ tcf}
(\prod_{i < \kappa}
\lambda_i / J^{bd}_\kappa)$, and \footnote{for the rest of this case
``$\lambda = \mu^+$" is not used; also $J^{bd}_\kappa$ can be replaced
by any larger ideal} let $\langle f_\alpha:\alpha < \lambda
\rangle$ exemplifying this.  Without loss of generality $\dsize \bigcup
_{j < i} \lambda_j < f_\alpha(i) < \lambda_i$.  Let $g:\kappa \times \mu
\times \kappa \times \mu \rightarrow \mu$ be one to one and onto, let
$f^\delta_\alpha = f_{\text{otp}(\alpha \cap C_\delta)}$ for
$\alpha \in C_\delta$ and let
$A^\delta_\alpha = \{g(i,f^\delta_\alpha(i),j,f^\delta_\alpha(j)):
i,j < \kappa\}$. \newline
\medskip

If $\delta = \sup(E_{\zeta(*)} \cap \text{ nacc}(C_\delta))$ 
and $E_{\zeta(*)} \cap \text{ nacc}(C_\delta) \ne \emptyset \text{ mod }
J_\delta$ then (as $J_\delta$ is $\lambda$-complete) choose $Y_\delta \in 
J_\delta$ such that for each $i < \kappa,\varepsilon < \lambda_i$ we have
\medskip
\roster
\item "{$(*)$}"  $(\exists \beta)[\beta \in E_{\zeta(*)} \cap 
\text{ nacc}(C_\delta) \and \beta \notin Y_\delta \and f^\delta_\beta(i) 
= \varepsilon] \Rightarrow$ \newline

$\quad \{\beta:\beta \in E_{\zeta(*)} \cap \text{ nacc}(C_\delta)
\and f^\delta_\beta(i) = \varepsilon \} \ne \emptyset \text{ mod }J_\delta$.
\endroster
\medskip

\noindent
Choose $i(\delta) < \kappa$ such that

$$
B^0_\delta =: \{ f^\delta_\beta(i(\delta)):\beta \in E_{\zeta(*)} \cap 
\text{ nacc}(C_\delta) \text{ and } \beta \notin Y_\delta\}
$$
\medskip

\noindent
is unbounded in $\lambda_i$.

Let $\xi_\varepsilon = \xi^\delta_\varepsilon$ be the $\varepsilon$-th 
member of $B^0_\delta$, for
$\varepsilon < \kappa$.  For each such $\varepsilon < \kappa$ for some 
$j_\varepsilon = j^\delta_\varepsilon \in (i(\delta) + 1 + \varepsilon,
\kappa)$ we have 
$B^{1,\delta}_\varepsilon =: \{ f^\delta_\beta(j_\varepsilon):
f^\delta_\beta(i(\delta)) = \xi^\delta_\varepsilon$ and \newline
$\beta \in
E_{\zeta(*)} \cap \text{ nacc}(C_\delta) \text{ and } \beta \notin Y_\delta\}$
is unbounded in  $\lambda_{j^\delta_\varepsilon}$.

Let $h_{\delta,\varepsilon}$ be a one to one function from
$[\dsize \bigcup_{j < \varepsilon} \lambda_j,\lambda_\varepsilon)$ into
$B^{1,\delta}_\varepsilon$.
\medskip

\noindent
Lastly we define $h_\delta$ as follows:

$$
\align
\text{\underbar{if} } \beta \in C_\delta,&\varepsilon < \kappa,
f^\delta_\beta(i(\delta)) = \xi^\delta_\varepsilon 
\text{ and } h_{\delta,\varepsilon}(\gamma)
= f^\delta_\beta(j^\delta_\varepsilon) \\
  &(\text{so } \gamma \in [\dsize \bigcup_{j < \varepsilon} \lambda_j,
\lambda_\varepsilon)) \text{ \underbar{then} } h_\delta(\beta) = \gamma
\endalign
$$
\smallskip

\noindent
and $h_\delta(\beta) = 0$ otherwise.  The rest is similar to the regular
case.  \hfill$\square_{\scite{3.1}}$
\bigskip

\proclaim{\stag{3.2} Claim}  If $\lambda = \mu^+,\mu$ regular uncountable and
$S \subseteq \{\delta  < \lambda:\text{cf}(\delta) = \mu \}$ is stationary
\underbar{then} for some strict 
$S$-club system $\bar C$ with $C_\delta = \{\alpha_{\delta,\zeta }:\zeta  < 
\mu \}$, (where $\alpha _{\delta,\zeta}$ is strictly increasing continuous 
in $\zeta)$ for every club  $E \subseteq \lambda$ for
stationarily many  $\delta \in S$,

$$
\{ \zeta < \mu:\alpha_{\delta,\zeta +1} \in  E \} \text{ is stationary (as 
a subset of } \mu).
$$
\endproclaim
\bigskip

\remark{\stag{3.2A} Remark}  1) If $S \in I[\lambda]$ then without 
loss of generality we can demand $(a)$ or we can demand $(b)$ (but 
not necessarily both), where  
\medskip
\roster
\item "{(a)}"  $X_\alpha = \{C_\delta \cap \alpha:\delta \in S,
\text{ is such that }
\alpha \in \text{ nacc}(C_\delta)\}$ has cardinality $\le \lambda$,
\item "{(b)}"  $\alpha \in \text{ nacc}(C_\delta) \Rightarrow C_\alpha =
C_\delta \cap \alpha$ but the conclusion is weakened to: \newline
for every club $E$ of $\lambda$ for stationarily many $\delta \in S$ 
the set \newline
$\{ \zeta < \mu:
(\alpha_{\delta,\zeta},\alpha_{\delta,\zeta + 1}) \cap E \ne \emptyset \}$ is
stationary.
\endroster
\medskip

\noindent
2)  In contrast to \cite[3.4]{Sh:413} here we allow $\mu$ inaccessible.\newline
3)  Clearly \scite{3.1}(2) can be applied to the results of \scite{3.2} 
i.e. with

$$
J_\delta = \biggl\{ A \subseteq C_\delta:\{\zeta < \lambda:\alpha
_{\delta,\zeta + 1} \notin A\} \text{ is not stationary}\biggr\}.
$$
\endremark
\bigskip

\demo{Proof}  We know that for some strict $S$-club system  $\bar C^0 = 
\langle C^0_\delta:\delta  \in  S \rangle$  we have  $\lambda \notin  
\text{ id}_p(\bar C^0)$ (see \cite[2.3(1)]{Sh:365}).  Let $C^0_\delta  = 
\{\alpha ^\delta _\zeta :\zeta  < \mu \}$  (increasing continuously in
$\zeta$).
We shall prove below that for some sequence of functions  
$\bar h = \langle h_\delta:\delta \in  S \rangle$, $h_\delta:\mu  
\rightarrow \mu$ we have 
\roster
\medskip
\item "{$(*)_{\bar h}$}"  for every club $E$ of $\mu^+$
for stationarily many  $\delta \in  S \cap \text{ acc}(E)$, \newline
the following subset of $\mu$ is stationary: \newline
$A^{\delta,*}_E =: \biggl\{
\zeta  < \mu:\alpha^\delta_\zeta \in  E \text{ and some ordinal in } 
\{\alpha^\delta_\xi:\zeta < \xi \le h_\delta(\zeta)\}$
\newline

$\qquad \qquad \qquad$ belongs to  $E \biggr\}$.
\endroster
\medskip

\noindent
The proof now breaks into two parts. \newline
\underbar{Proving} $(*)_{\bar h}$ \underbar{suffices}\,.
\medskip
\noindent
For each  club $E$ of $\lambda$, let
$Z_E =: \{ \delta \in S:\delta = \sup(E \cap \text{ nacc}(C^0_\delta))\}$,
and note that this
set is a stationary subset of $\lambda$ (by the choice of $\bar C^0$).
For each such $E$ and $\delta \in Z_E$ let $f_{\delta,E}$ be 
the partial function from $\mu$ to $\mu$ defined by

$$
f_{\delta,E}(\zeta) = \text{ Sup}\{\xi:\zeta < \xi \le h_\delta(\zeta) 
\text{ and } \alpha^\delta_\xi \in E\}.
$$
\medskip

\noindent
So if there is no such $\xi$, then $f_{\alpha,E}(\zeta)$ is not well
defined (i.e. if the supremum is on the empty set) but if $\xi =
f_{\alpha,E}(\zeta)$ is well defined then $\alpha^\delta_\xi \in E,
\xi \le h_\delta(\zeta)$ (because $\alpha^\delta_\xi$ is increasing
continuous in $\xi$ and $E$ is a club of $\lambda$).  Let \newline
$Y_E =: \{ \delta \in Z_E:\text{Dom}(f_{\delta,E})$ is a
stationary subset of $\mu\}$.  So by $(*)_{\bar h}$, we know that
\medskip
\roster
\item "{$\bigoplus$}"  for every club $E$ of $\mu^+$ the set $Y_E$ is a
stationary subset of $\mu^+$.
\endroster
\medskip

\noindent
Also
\medskip
\roster
\item "{$\bigotimes_1$}"  if $E_2 \subseteq E_1$ are clubs of $\mu^+$ then
$Z_{E_2} \subseteq Z_{E_1}$ and $Y_{E_2} \subseteq Y_{E_1}$ and for \newline
$\delta \in Y_{E_2},\text{Dom}(f_{\delta,E_2}) \subseteq \text{ Dom}(
f_{\delta,E_1})$ and \newline
$\zeta \in \text{ Dom}(f_{\delta,E_2}) \Rightarrow
f_{\delta,E_2}(\zeta) \le f_{\delta,E_1}(\zeta)$.
\endroster
\medskip

\noindent
We claim that
\medskip
\roster
\item "{$\bigotimes_2$}"  for some club $E_0$ of $\mu^+$ for every club
$E \subseteq E_0$ of $\mu^+$ for stationarily many $\delta \in S$ we
have
{\roster
\itemitem{ (i) }  $\delta = \text{ sup}(E \cap \text{ nacc } C_\delta)$,
\itemitem{ (ii) }  $\{ \zeta < \mu:\zeta \in \text{ Dom}(f_{E,\delta})$
(hence $\zeta \in \text{ Dom }f_{E_0,\delta})$ and \newline

$\qquad \qquad \quad f_{E,\delta}(\zeta) = 
f_{E_0,\delta}(\zeta)\}$ is a stationary subset of $\mu$.
\endroster}
\endroster
\medskip

\noindent
If this fails, then for any club $E_0$ of $\lambda$ there is a club
$E(E_0) \subseteq E_0$ of $\lambda$, such that

$$
\align
A_{E_0} = \biggl\{ \delta:&\delta \in S,\delta = \text{ sup}(E(E_0) \cap
\text{ nacc}(C_\delta)) \text{ and for some club} \\
  &e_{E_0,\delta} \text{ of } \mu \text{ we have}\\
  &\zeta \in e_{E_0,\delta} \cap \text{ Dom}(f_{E(E_0),\delta})
\Rightarrow f_{E(E_0),\delta}(\zeta) = f_{E_0,\delta}(\zeta) \biggr\}
\endalign
$$
\medskip

\noindent
is not a stationary subset of $\lambda = \mu^+$.  By obvious 
monotonicity we can replace
$E(E_0)$ by any club of $\mu^+$ which is a subset of it, so without loss of
generality $A_{E_0} = \emptyset$.

By induction on $n < \omega$ choose clubs $E_n$ of $\mu^+$ such that
$E_0 = \mu^+$ and \newline
$E_{n+1} = E(E_n)$. \newline
Then $E_\omega =: \dsize \bigcap_{n < \omega} E_n$ is a club of $\mu^+$ and,
by $\bigoplus$ above, $Y_{E_\omega} \subseteq S$ is a stationary subset of
$\lambda$, so we can choose a $\delta(*) \in Y_{E_\omega}$.  So
$f_{E_\omega,\delta(*)}$ has domain a stationary subset of $\mu$ (see the
definition of $Y_{E_\omega})$ and by $\bigotimes_1$ we know that \newline
$n < \omega \Rightarrow \text{ Dom}(f_{E_\omega,\delta(*)}) \subseteq
\text{ Dom}(f_{E_n,\delta(*)})$.  Also there is an $e_{E_n,\delta(*)}$, a club
of $\mu$, such that

$$
\zeta \in e_{E_n,\delta(*)} \cap \text{ Dom}(f_{E_{n+1},
\delta(*)}) \Rightarrow
f_{E_{n+1},\delta(*)}(\zeta) < f_{E_n,\delta(*)}(\zeta)
$$
\medskip

\noindent
(see the choice of $E_{n+1} = E(E_n)$ i.e. the function $E$).  So
$e_{\delta(*)} =: \dsize \bigcap_{n < \omega} e_{E_n,\delta(*)}$ is a club 
of $\mu$ and,
as $\text{Dom}(f_{E_\omega,\delta(*)})$ is a stationary subset of $\mu$, we
can find \newline
$\zeta(*) \in e_{\delta(*)} \cap \text{ Dom}(f_{E_\omega,\delta(*)})$,
hence 
$\zeta(*) \in \dsize \bigcap_{n < \omega} \text{ Dom}
(f_{E_n,\delta(*)}) \cap \dsize \bigcap_{n < \omega} e_{E_n,\delta(*)}$,
so that
$\langle f_{E_n,\delta(*)}(\zeta(*)):n < \omega \rangle$ is a well defined
strictly increasing $\omega$-sequence of ordinals - a contradiction.
So $\bigotimes_2$ cannot fail, and this gives the desired conclusion. 
\enddemo
\bigskip

\demo{\underbar{Proof of} $(*)_{\bar h}$ \underbar{holds for some} $\bar h$}
\newline

So assume that for no $\bar h$  does $(*)_{\bar h}$ holds, hence (by shrinking
$E$) we can assume that for every $\bar h = \langle h_\delta:\delta \in S
\rangle,h_\delta:\mu \rightarrow \mu$, for some club $E$ for \underbar{every}
$\delta \in S,A^{\delta,*}_E$ is not stationary (in $\mu$).  By induction on
$n < \omega$, we define $E_n$, \newline
$\bar h^n = \langle h^n_\delta:
\delta \in S \rangle, \bar e^n = \langle e^n_\delta:\delta \in S \rangle$, 
with $E_n$ a club of $\lambda,e^n_\delta$ club of $\mu,h^n_\delta:\mu 
\rightarrow \mu$ as follows. \newline
\noindent
Let  $E_0 = \lambda$,  $h^0_\delta(\zeta) = \zeta + 1$ and  
$e^n_\delta  = \mu$. \newline
\noindent
If  $E_0,...,E_n$,  $\bar h^0,...,\bar h^n$,  $\bar e^0,...,\bar e^n$ are 
defined, necessarily $(\ast)_{\bar h^n}$ fail, so for some club $E_{n+1}$ 
of $\lambda$  for every  $\delta  \in  S \cap \text{ acc}(E_{n+1})$
there is a club  
$e^{n+1}_\delta \subseteq \text{ acc}(e^n_\delta)$ of $\mu$, such that

$$
\zeta \in e^{n+1}_\delta \Rightarrow \{ \alpha^\delta_\xi:\zeta < \xi \le
h_\delta(\zeta) \} \cap E_{n+1} = \emptyset.
$$
\medskip

\noindent
Choose  $h^{n+1}_\delta:\mu \rightarrow \mu$  such that  
$(\forall \zeta < \mu)(h^n_\delta(\zeta) < h^{n+1}_\delta(\zeta))$ and
if \newline
$\delta = \sup(E_{n+1} \cap \text{ nacc}(C_\delta))$ then
$\zeta < \mu \Rightarrow \{ \alpha^\delta_\xi:\zeta < \xi \le h^{n+1}
_\delta(\zeta)\} \cap E_{n+1} \ne \emptyset$.
\medskip

\noindent
There is no problem to carry out this inductive definition.  
By the choice of $\bar C^0$, for some
$\delta \in \text{ acc}(\dsize \bigcap_{n<\omega} E_n)$, we have $\delta =
\sup(A')$, where \newline
$A^\prime =: (\text{acc} \dsize \bigcap_{n < \omega}
 E_n) \cap \text{ nacc}(C^0_\delta)$.  Let $A \subseteq \mu$ be such that
$A^\prime = \{\alpha^\delta_\zeta:\zeta \in A\}$ (remember
$\alpha^\delta_\zeta$ is increasing with $\zeta$) and
let $\zeta$ be the second member of $\dsize \bigcap_{n < \omega} e^n_\delta$.
As $A'$ is unbounded in $\delta$, clearly $A$ is unbounded in $\mu$ and
$\dsize \bigcap_{n < \omega} e^n_\delta$ is a club of $\mu$ as
$\mu = \text{cf}(\mu) > \aleph_0$.  Also as $A' \subseteq \text{ nacc}
(C^0_\delta)$ clearly $A$ is a set of successor ordinals (or zero).

Note that $\sup(e^\delta_n \cap \zeta)$ is well defined (as
$\text{Min}(e^\delta_n) \le \text{ Min}(\dsize \bigcap_{n < \omega} 
e^n_\delta) < \zeta)$ and $\sup(e^\delta_n \cap
\zeta) < \zeta$ (as $\zeta$ is a successor ordinal).
Now  $\langle \sup (e^\delta_n \cap  
\zeta ):n < \omega \rangle$  is non-increasing (as  $e_\delta^n$ decreases 
with $n$), hence for some $n(*) < \omega$ we have $n > n(*) \Rightarrow  
\sup (e_\delta^n \cap \zeta) = \sup(e_\delta^{n(\ast)} \cap \zeta$) and  
call this ordinal $\xi$ so that $\xi \in e^\delta_{n(*)+1}$ and
$h^{n(*)}_\delta(\xi) = h^{n(*)+1}_\delta(\xi)$, so we get a contradiction
for $n(*) + 1$. \newline
So $(*)_{\bar h}$ holds for some $\bar h$, which suffices, as indicated
above. \hfill$\square_{\scite{3.2}}$
\enddemo
\bigskip

\demo{\stag{3.2B} Discussion}  1)  We can squeeze a little more, 
but it is not so clear if with much gain.  So assume
\medskip
\roster
\item "{$(*)_0$}"  $\mu$ is regular uncountable, $\lambda = \mu^+,
S \subseteq \{ \delta < \lambda:\text{cf}(\delta) = \mu \}$ stationary, $I$
an ideal on $S$, $\bar C = \langle C_\delta:\delta \in S \rangle$ a strict
$S$-club system, $\bar J = \langle J_\delta:\delta \in S \rangle$ with
$J_\delta$ an ideal on $C_\delta$ extending $J^{bd}_{C_\delta} +
(\text{acc}(C_\delta))$, such that for any club $E$ of $\lambda$ we have
$\{ \delta
\in S:E \cap C_\delta \ne \emptyset \text{ mod } J_\delta\} \ne \emptyset 
\text{ mod }I$.
\endroster
\medskip

\noindent
2)  If we immitate the proof of \scite{3.2} we get
\medskip
\roster
\item "{$(*)_1$}"  if for $\delta \in S,J_\delta$ is not $\chi$-regular
(see the definition below) and $\chi \le \mu$ then we can find 
$\bar e = \langle e_\delta:
\delta \in S \rangle$ and $\bar g = \langle g_\delta:\delta \in S \rangle$ 
such that
\item "{$(*)'_1$}"  $e_\delta$ is a club of $\delta,e_\delta \subseteq
\text{ acc}(C_\delta),g_\delta:\text{nacc}(C_\delta) \backslash (\text{min}
(e_\delta) + 1) \rightarrow e_\delta$ is defined by $g_\delta(\alpha) = \sup
(e_\delta \cap \alpha)$ and for every club $E$ of $\lambda$
$$
\align
\biggl\{ \delta \in S:&\,E \cap \text{ nacc}(C_\delta) \ne \emptyset
\text{ mod } J_\delta \text{ and} \\
  &\text{ Rang}(g_\delta \restriction (E \cap \text{ nacc}(C_\delta)))
\text{ is a stationary subset of } \delta \biggr\} \ne \emptyset \text{ mod }
I.
\endalign
$$
\endroster
\medskip

\noindent
3) \underbar{Definition}: An ideal $J$ on a set $C$ is $\chi$-regular if
there is a set $A \subseteq C$, \newline
$A \ne \emptyset \text{ mod }J$ and a function
$f:A \rightarrow [\chi]^{< \aleph_0}$ such that \newline
$\gamma < \chi \Rightarrow
\{ x \in A:\gamma \notin f(x)\} = \emptyset \text{ mod } J$. \newline
If $\chi = |C|$, we may omit it. \newline
[How do we prove $(*)'_1$?  Try $\chi$ times $E_\zeta,\langle e^\zeta_\delta:
\delta \in S \rangle$ (for $\zeta < \chi$)]. \newline
4)  We can try to get results like \scite{3.1}.  Now
\medskip
\roster
\item "{$(*)_2$}"  assume $\lambda,\mu,S,I,\bar C,\bar J$ 
are as in $(*)_0$ and
$\bar e,\bar g$ as in $(*)'_1$ and $\kappa < \mu$ and for $\delta \in S,
J^0_\delta =: \{ a \subseteq e_\delta:\{ \alpha \in \text{ Dom}(g_\delta):
g(\alpha) \in a\} \in J_\delta\}$ is weakly normal and 
$\mu$ satisfies the condition from \cite[Lemma 2.12]{Sh:365}. \underbar{Then} 
we can find $h_\delta:e_\delta \rightarrow \kappa$ such that for every 
club $E$ of $\lambda$, \newline
$\{ \delta \in
S:\text{for each } \gamma < \kappa \text{ the set }\{ \alpha \in \text{ nacc}
(C_\delta):h_\delta(g_\delta(\alpha)) = \gamma\} \text{ is }$ \newline
$\ne \emptyset \text{ mod }J_\delta\} \ne \emptyset \text{ mod } I$.
\endroster
\medskip

\noindent
[Why?  For each $\delta \in S,\alpha \in \text{ acc}(e_\delta)$ choose
a club \newline
$d_{\delta,\alpha} \subseteq e_\delta \cap \alpha$ such that for no club
$d \subseteq e_\delta$ of $\delta$ do we \newline
have
$(\forall \gamma < \delta)(\exists \alpha \in \text{ acc}(e_\delta))[d \cap
\gamma \subseteq d_{\delta,\alpha}]$.  Now for every club $E$ of $\lambda$
let \newline
$S_E = \{ \delta:E \cap \text{ nacc}(C_\delta) \ne \emptyset \text{ mod }
J_\delta,\text{ and } g''_\delta(E \cap \text{ nacc}(C_\delta)) \text{ is
stationary}\}$ and for
$\delta \in E$ and $\varepsilon < \mu$, we choose by
induction on $\zeta < \kappa,\xi(\delta,\varepsilon)$ as the first
$\xi \in e_\delta$ such that: $\xi > \dsize \bigcup_{\zeta < \varepsilon} 
\xi(\delta,\zeta)$ and $\{ \alpha \in \text{ Dom}(g_\delta):
\alpha \in E$ and the $\varepsilon$-th member of $d_{\delta,g_\delta(\alpha)}$
is in the interval $[\dsize \bigcup_{\zeta < \varepsilon} 
\xi(\delta,\zeta),\xi)]\} \ne \emptyset \text{ mod } J_\delta$. \newline
5)  We deal below with successor of singulars and with inaccessibles, we
can do parallel things.
\enddemo
\bigskip

\proclaim{\stag{3.3} Claim}  Suppose 
$\mu$ is a singular cardinal of cofinality $\kappa,\kappa > \aleph_0$ 
and \newline
$S \subseteq \{ \delta < \mu^+:\text{cf}(\delta) = \kappa\}$ is
stationary, and $\bar C = \langle C_\delta:\delta \in S \rangle$ is an
$S$-club system satisfying $\mu^+ \notin \text{ id}^p(\bar C,\bar 
J^{b[\mu]})$ where $\bar J^{b[\mu]} = \langle J^{b[\mu]}_{C_\delta}:\delta \in
S \rangle$ and \newline
$J^{b[\mu]}_{C_\delta} =: \{ A \subseteq C_\delta:
\text{for some } \theta < \mu,\text{ we have } \delta > \sup\{\alpha \in A:
\text{cf}(\alpha) > \theta\}\}$. \underbar{Then} we can find a strict 
$\lambda$-club system
$\bar e^* = \langle e^*_\delta:\delta < \lambda \rangle$ such that
\medskip
\roster
\item "{$(*)$}"  for every club $E$ of $\mu^+$, for stationarily many
$\delta \in S$, for every $\alpha < \delta$ and $\theta < \mu$ for some
$\beta$ we have
{\roster
\itemitem{ $(**)_{E,\beta}$ }  $\beta \in \text{ nacc}(C_\delta)$ and
$\beta > \alpha$ and $\text{cf}(\beta) > \theta$ and \newline
$\{ \gamma \in e^*_\beta:\gamma \in E \text{ and min}(e^*_\beta \backslash
(\gamma + 1)) \text{ belongs to } E\}$ \newline
is a stationary subset of $\beta$.
\endroster}
\endroster
\endproclaim
\bigskip

\remark{\stag{3.3A} Remark}  1) We know 
that for the given $\mu$ and $S$ there is $\bar C$ as in 
the assumption by \cite[\S2]{Sh:365}.  Moreover, if
$\kappa > \aleph_0$ then there is such nice strict $\bar C$. \newline
2) Remember $J^{b[\mu]}_\delta = \{A \subseteq C_\delta:\text{for some }
\theta < \mu \text{ and } \alpha < \delta$ we have \newline

$\qquad \quad \qquad \qquad \qquad \qquad \qquad(\forall \beta
\in C_\delta)(\beta < \alpha \vee \text{ cf}(\beta) < \theta \vee \beta \in
\text{ nacc}(C_\delta))\}$.
\endremark
\bigskip

\demo{Proof} Let $\bar e = \langle e_\beta:\beta < \lambda \rangle$ be a
strict $\lambda$-club system where $e_\beta = \{ \alpha^\beta_\zeta:\zeta
< \text{ cf}(\beta)\}$ is a (strictly) increasing and continuous 
enumeration of $e_\beta$ (with limit $\delta$).  Now we claim that for some
$\bar h = \langle \bar h_\beta:\beta < \lambda,\beta \text{ limit}
\rangle$ with $h_\beta$ a function from $e_\beta$ to $e_\beta$ and $\dsize
\bigwedge_{\alpha \in e_\beta} h_\beta(\alpha) > \alpha$, we have
\medskip
\roster
\item "{$(*)_{\bar h}$}"  for every club $E$ of $\mu^+$, for stationarily
many $\delta \in S \cap \text{ acc}(E),A^\delta_E \notin J^{b[\mu]}
_{C_\delta}$ where $A^\delta_E$ is the set of all $\beta \in C_\delta$ such
that the following subset of $e_\beta$ is stationary (in $\beta$):
$$
\{\gamma \in e_\beta:\gamma \in E \text{ and min}(e_\beta \backslash
(\gamma + 1)) \in E \}.
$$
\endroster
\medskip

\noindent
The rest is like the proof of \scite{3.2} repeating $\kappa^+$ times instead
$\omega$ and using ``$J^{b[\mu]}_{C_\delta}$ is $(\le \kappa)$-based".
\hfill$\square_{\scite{3.3}}$
\enddemo
\bigskip

\proclaim{\stag{3.4} Claim}  Suppose $\lambda$ is inaccessible, $S \subseteq
\lambda$ is a stationary set of inaccessibles, $\bar C$ an $S$-club system
such that $\lambda \notin \text{ id}^p(\bar C)$.  \underbar{Then} we can find
$\bar h = \langle h_\delta:\delta \in S \rangle$ with $h_\delta:C_\delta
\rightarrow C_\delta$, such that $\alpha < h(\alpha)$ and
\medskip
\roster
\item "{$(*)$}"  for every club $E$ of $\lambda$, for stationarily many
$\delta \in S \cap \text{ acc}(E)$ we have that
$$
\{ \alpha \in C_\delta:\alpha \in E \text{ and } h(\alpha) \in E\}
\text{ is a stationary subset of } \delta.
$$
\endroster
\medskip

\noindent
So for some $C'_\delta = \{ \alpha_{\delta,\zeta}:\zeta < \delta\} \subseteq
C_\delta,\alpha_{\delta,\zeta}$ increasing continuous in $\zeta$ we have
$h(\alpha_{\delta,\zeta}) = \alpha_{\delta,\zeta + 1}$. 
\endproclaim
\bigskip

\remark{Remark}  Under quite mild conditions on $\lambda$ and $S$ there is
$\bar C$ as required - see \newline
\cite[2.12,p.134]{Sh:365}.
\endremark
\bigskip

\demo{Proof}  Like \scite{3.2}.
\enddemo
\bigskip

\proclaim{\stag{3.5} Claim}  Let $\lambda = \text{ cf}(\lambda) > \aleph_0,
S \subseteq \lambda$ stationary, 
$D$ a normal $\lambda^+$-saturated filter on $\lambda$, $S$ is 
$D$-positive (i.e. $S \in D^+$, $\lambda \backslash S \notin D$). \newline
1)  Assume $\langle(C_\delta,I_\delta):\delta \in S \rangle$ is such that
\medskip
\roster
\item "{(a)}"  $C_\delta \subseteq \delta = \sup(C_\delta),I_\delta
\subseteq {\Cal P}(C_\delta)$,
\item "{(b)}"  for every club $E$ of $\lambda$,
$$
\{ \delta \in S:\text{for some } A \in I_\delta \text{ we have }
\delta > \sup(A \backslash E)\} \in D^+.
$$
Then for some stationary $S_0 \subseteq S,S_0 \in D^+$ we have
\item "{$(b)^+$}"  for every club $E$ of $\lambda$
$$
\{ \delta \in S:\text{for no } A \in I \text{ do we have }
\delta > \sup(A \backslash E)\} = \emptyset \text{ mod } D.
$$
\endroster
\medskip

\noindent
2)  Assume $\langle {\Cal P}_\delta:\delta \in S \rangle$ is such that (here
really presaturated is enough)
\medskip
\roster
\item "{$(*)$}"  for every $D$-positive $S_0 \subseteq S$ for some
$D$-positive $S_1 \subseteq S_0$ and \newline
$\langle (C_\delta,I_\delta):\delta \in S \rangle$ we have
$(C_\delta,I_\delta) \in {\Cal P}_\delta,C_\delta \subseteq
\delta = \sup(C_\delta),I_\delta \subseteq {\Cal P}(C_\delta)$ and for
every club $E$ of $\lambda$ \newline
$\{ \delta \in S_1:\text{for some } A \in I_\delta,\delta > \sup(A \backslash
E)\} \ne \emptyset \text{ mod } D$.
\endroster
\medskip

\noindent
Then
\medskip
\roster
\item "{$(**)$}"  for some $\langle (C_\delta,A_\delta):\delta \in
S \rangle$ we have $(C_\delta,I_\delta) \in {\Cal P}_\delta,
C_\delta \subseteq \delta = \sup(C_\delta)$, \newline
$I_\delta \subseteq {\Cal P}(C_\delta)$ and for every club $E$ of $\lambda$
$$
\{ \delta \in S:\text{ for no } A \in I_\delta,\delta > \sup(A \backslash
E)\} = \emptyset \text{ mod } D.
$$
\endroster
\endproclaim
\bigskip

\remark{Remark}  This is a straightforward 
generalization of \cite[III,\S6.2B]{Sh:e}.  Independently Gitik found 
related results on generic extensions which
were continued in \newline
\cite{DjSh:562} and in \cite{GiSh:577}.
\endremark
\bigskip

\demo{Proof}  The same.
\enddemo
\bigskip

\proclaim{\stag{3.6} Lemma}  Suppose $\lambda$ is regular uncountable and
$S \subseteq \{ \delta < \lambda^+:\text{cf}(\delta) = \lambda\}$ is
stationary.  \underbar{Then} we can find $\langle (C_\delta,h_\delta,\chi_
\delta):\delta \in S \rangle$ and $D$ such that
\medskip
\roster
\item "{$(A)$}"  $D$ is a normal filter on $\lambda^+$,
\item "{$(B)$}"  $C_\delta$ is a club of $\delta$, say $C_\delta = \{
\alpha_{\delta,\zeta}:\zeta < \lambda\}$, with $\alpha_{\delta,\zeta}$
increasing continuous in $\zeta$,
\item "{$(C)$}"  $h_\delta$ is a function from $C_\delta$ to $\chi_\delta,
\chi_\delta \le \lambda$,
\item "{$(D)$}"  if $A \in D^+$ (i.e. $A \subseteq \lambda^+ \and \lambda^+
\backslash A \notin D$) and $E$ is a club of $\lambda^+$, then the following
set belongs to $D^+$:
$$
\align
B_{E,A} =: \biggl\{ \delta:&\,\delta \in A \cap S,\delta \in \text{ acc}(E)
\text{ and for each } i < \chi_\delta \\
  &\,\{\zeta < \lambda:\alpha_{\delta,\zeta +1} \in E \text{ and }
h_\delta(\alpha_{\delta,\zeta}) = i \\
  &\,\text{(and } \alpha_{\delta,\zeta} \in E)\} \text{ is a stationary subset
of } \lambda \biggr\}
\endalign
$$
\noindent
(hence, for some $\alpha < \lambda^+$ and $\zeta < \lambda$, the set \newline
$B_{E,A,\alpha} =: \{ \delta \in B_{E,A}:\alpha = \alpha_{\delta,\zeta}\}$
is in $D^+$).
\item "{$(E)$}"  If $\gamma < \lambda^+$ and $\chi$ satisfies one of the
conditions listed below, then \newline
$S_{\gamma,\chi} = \{ \delta \in
S:\gamma = \text{ Min}(C_\delta)$ and $\chi_\delta = \chi\} \in D^+$ where
{\roster
\itemitem{ $(\alpha)$ }  $\lambda = \chi^+$,
\itemitem{ $(\beta)$ }  $\lambda$ is inaccessible not strongly inaccessible,
$\chi < \lambda$ and there is $T$ such that \newline

$\qquad (a) \quad T$ is a tree with $< \lambda$ nodes and a set $\Gamma$ of
branches, $|\Gamma| = \lambda$, \newline

$\qquad (b)'$ \,\,if $T' \subseteq T,T'$ downward closed 
and $(\exists^\lambda \eta \in \Gamma)$ \newline

$\qquad \qquad (\eta$ a branch of $T'$) 
\underbar{then} $T'$ has an antichain of cardinality $\ge \chi$,
\newline
\itemitem{ $(\gamma)$ }  $\lambda$ is inaccessible not 
strongly inaccessible and \newline
$\chi = \text{ Min}\{\chi:\text{for some } \theta \le \chi
\text{ we have } \chi^\theta \ge \lambda\}$,
\itemitem{ $(\delta)$ }  $\lambda$ is strongly inaccessible not 
ineffable; i.e. $\lambda$ is Mahlo and \newline
we can find $\bar A = \langle A_\mu:\mu < \lambda$ is inaccessible $\rangle$,
\newline
$A_\mu \subseteq \mu$ so that for no stationary
$\Gamma \subseteq \{ \mu < \lambda:\mu \text{ inaccessible}\}$ \newline
and $A \subseteq \lambda$ do we have:
$\mu \in \Gamma \Rightarrow A_\mu = A \cap \mu$.
\endroster}
\endroster
\endproclaim
\bigskip

\remark{\stag{3.6A} Remark}  We can replace $\lambda^+$ in \scite{3.6} and
any $\mu = \text{ cf}(\mu) > \lambda$, as if $\mu > \lambda^+$ we have even
a stronger theorem.
\endremark
\bigskip

\demo{Proof}  Let for $\lambda = \text{ cf}(\lambda) > \aleph_0$,
$$
\align
\Theta = \Theta_\lambda = \biggl\{
\chi \le \lambda:&\text{ if } S' \subseteq \{ \delta < \lambda^+:
\text{cf}(\delta) = \lambda\} \text{ is stationary} \\
  &\text{ \underbar{then} we can find } \langle (C_\delta,h_\delta):\delta
\in S' \rangle \text{ such that} \\
  &\,(a) \quad C_\delta \text{ is a club of } \delta \text{ of order type }
\lambda, \\
  &\,(b) \quad h_\delta:C_\delta \rightarrow \chi, \\
  &\,(c) \quad \text{for every club } E \text{ of } \lambda^+ \text{ for
stationarily many} \\
  &\qquad \delta \in S' \cap \text{ acc}(E) \text{ we have:} \\
  &\qquad i < \chi \Rightarrow B_E = \{\alpha \in C_\delta:\alpha \in E,
h(\alpha) = i \text{ and} \\
  &\qquad \qquad \qquad \qquad \qquad \qquad \quad \text{ min}
(C_\delta \backslash (\alpha + 1)) \in E \} \\
  &\qquad \text{ is a stationary subset of } \delta \biggr\}.
\endalign
$$
\medskip

\noindent

\medskip

Now we first show
\roster
\item "{$\bigotimes$}"  for each of the cases from clause (E), the $\chi$
belongs to $\Theta$.
\endroster
\enddemo
\bigskip

\demo{Proof of sufficiency of $\bigotimes$}  We can partition $S$ to
$\lambda^+$ stationary sets so we can find a partition $\langle
S_{\chi,\alpha}:\chi \in \Theta \text{ and } \alpha < \lambda^+ \rangle$ of
$S$ to stationary sets.  Without loss of generality, $\alpha \le
\text{ Min}(S_{\chi,\alpha})$ and let $\langle (C^0_\delta,h^0_\delta):\delta
\in S_{\chi,\alpha} \rangle$ be as guaranteed by ``$\chi \in \Theta$" for
the stationary set $S_{\chi,\alpha}$.  Now define $C_\delta,h_\delta$ for
$\delta \in S$ by:

$C_\delta$ is $C^0_\delta \cup\{ \alpha\} \backslash \alpha$
\underbar{if} $\delta \in S_{\chi,\alpha}$ and $\alpha < \delta,
h_\delta(\beta)$ is $h_\delta(\beta)$ if $\beta \in C_\delta \cap
C^0_\delta$ and is zero otherwise.  Of course, $\chi_\delta = \chi$ if
$\delta \in S_{\chi,\alpha}$.
\medskip

\noindent
Lastly, let
$$
\align
D = \biggl\{ A \subseteq \lambda^+:&\text{ for some club } E \text{ of }
\lambda^+, \text{ for every} \\
  &\,\delta \in S \cap \text{ acc}(E) \backslash A \text{ for some } i <
\chi_\delta, \\
  &\text{ the set } \{ \beta \in C_\delta:\beta \in E,h_\delta(\beta) = i
\text{ and \, min}(C_\delta \backslash (\beta + 1) \in E \} \\
  &\,\text{is not a stationary subset of } \delta \biggr\}.
\endalign
$$
\medskip

\noindent
So $D$ and $\langle(C_\delta,h_\delta,\chi_\delta):\delta \in S \rangle$
have been defined, and we have to check clauses (A)-(E). \newline
Note that $\Theta \ne \emptyset$ and the proof which appears later does not
rely on the intermediate proofs.
\enddemo
\bigskip

\noindent
\underbar{Clause $(A)$}:  Suppose $A_\zeta \in D$ for $\zeta < \lambda$, so
for each $\zeta$ there is a club $E_\zeta$ of $\lambda^+$
\medskip
\roster
\item "{$(*)$}"  if $\delta \in S_{\chi,\gamma}$ and
$\delta \in S \cap \text{ acc}(E) \backslash A_\zeta$ then \newline
$\{ \alpha \in C_\delta:\alpha \in E,\text{Min}(C_\delta \backslash
(\alpha + 1)) \in E \text{ and } h_\delta(\alpha) = i_\zeta\}$
is not stationary in $\delta$.
\endroster
\medskip

\noindent
Clearly clubs of $\lambda^+$ belong to $D$. \newline
Clearly $A \supseteq A_\zeta \Rightarrow A \in D$ (by the definition),
witnessed by the same $E_\zeta$. \newline
Also $A = A_0 \cap A_1 \in D$ as witnessed by $E = E_0 \cap E_1$. \newline
Lastly, $A = \underset{\zeta < \lambda} {}\to \triangle \,\,A_\zeta = 
\{ \alpha < \lambda^+:
\alpha \in \dsize \bigcap_{\zeta < 1 + \alpha} A_\zeta\}$ belong to $D$ as
witnessed by $E = \{ \alpha < \lambda^+:\alpha \in \dsize \bigcap_{\zeta <
1 + \alpha} E_\zeta\}$.  Note that if $\delta \in S \cap \text{ acc}(E)
\backslash A$ then for some $\zeta < \delta$ 

$$
\delta \in S \cap \text{ acc}(E) \backslash A_\zeta \subseteq (S \cap
\text{ acc}(E_\zeta) \backslash A_\zeta) \cup (1 + \zeta)
$$
\medskip

\noindent
as $E_\zeta \backslash E$ is a bounded subset of $\delta$; included in
$1 + \zeta$ so from the
conclusion of $(*)$ for $\delta,A_\zeta,E_\zeta$ we get it for $\zeta,A,E$.
\smallskip

Lastly $\emptyset \notin D$; otherwise, let $E$ be a club of $\lambda^+$ 
witnessing it, i.e. $(*)$ holds in this case.  Choose $\chi \in \Theta$ and
$\alpha = 0$ and use on it the choice of $\langle C^0_\delta:\delta \in
S_{\chi,0} \rangle$ to show that for some $\delta \in S_{\chi,0} \subseteq
S$ contradict the implication in $(*)$.
\bigskip

\noindent
\underbar{Clause $(B)$}:  Trivial.
\bigskip

\noindent
\underbar{Clause $(C)$}:  Trivial.
\bigskip

\noindent
\underbar{Clause $(D)$}:  Note that we can ignore the ``$\alpha_{\delta,
\zeta} \in E$" as $\delta \in \text{ acc}(E)$ implies that it holds for a
club of $\zeta$'s.
Assume $A \in D^+$ (for clause $(A)$) and $E$ is a
club of $\lambda^+$, which contradicts clause $(D)$ so $B_{E,A} \notin D^+$,
hence $\lambda^+ \backslash B_{E,A} \in D$.  Also $E$ witnessed that
$\lambda^+ \backslash (A \backslash B_{E,A}) \in D$ by the definition of $D$.
But by clause $(A)$ we know $D$ is a filter on $\lambda^+$ so $(\lambda^+
\backslash B_{E,A}) \cap (\lambda^+ \backslash (A \backslash B_{E,A})$ belong
to $D$, but this is the set $\lambda^+ \backslash B_{E,A} \backslash 
(A \backslash B_{E,A})$ which is (as $B_{E,A} \subseteq A$ by its 
definition) just $\lambda \backslash A$.  
So $\lambda \backslash A \in D$ hence $A \notin D^+$ - a contradiction.
\bigskip

\noindent
\underbar{Clause $(E)$}:  By the proof of $\emptyset \notin D$ above, if
$\chi \in \Theta$, also $S_{\chi,\alpha} \in D^+$, and by the definition of
$\bar C,\bar C \restriction S_{\chi,\alpha}$ is as required.  So it is enough
to show
\bigskip

\proclaim{\stag{3.7} Claim}  If 
$\chi < \lambda = \text{ cf}(\lambda)$ and $\chi$
satisfies one of the clauses of \scite{3.6(E)}, then 
$\chi \in \Theta$ (from the proof of \scite{3.6}).
\endproclaim
\bigskip

\demo{Proof}  \newline
\underbar{Case $(\alpha)$}:  By \scite{3.1}.
\bigskip

\noindent
\underbar{Case $(\beta)$}:  Like the proof of \scite{3.1}, for 
more details see \cite[\S3]{Sh:413}.
\bigskip

\noindent
\underbar{Case $(\gamma)$}:  This is a particular case of case $(\beta)$.
\bigskip

\noindent
\underbar{Case $(\delta)$}:  Similar proof (or use \scite{3.8}).
\hfill$\square_{\scite{3.7}}\square_{\scite{3.6}}$
\enddemo
\bigskip

\noindent
More generally (see \cite{Sh:413}):
\proclaim{\stag{3.8} Claim}  Let $\lambda = \text{ cf}(\lambda) > \chi$.  A
sufficient condition for $\chi \in \Theta_\lambda$ is the existence of some
$\zeta < \lambda^+$ such that
\medskip
\roster
\item "{$\bigotimes$}"  in the following game of length $\zeta$, first
player has no winning strategy: \newline
in the $\varepsilon$-th move first player chooses a function $f_\varepsilon:
\lambda \rightarrow \chi$ and second player chooses $\beta_\varepsilon <
\chi$.  In the end, first player wins the play if \newline
$\{ \alpha < \lambda:
\text{for every } \varepsilon < \gamma,f_\varepsilon(\alpha) \ne 
\beta_\varepsilon\}$ is a stationary subset of $\lambda$.
\endroster
\medskip

\noindent
(If we weaken the demand in $\Theta_\lambda$ from stationary to unbounded
in $\lambda$, we can weaken it here too).
\endproclaim
\newpage

\head {\S4 More on $Pr_6$} \endhead
\resetall
\bigskip

\proclaim{\stag{4.1} Claim}  $Pr_6(\lambda^+,\lambda^+,\lambda^+,\lambda)$ for
$\lambda$ regular.
\endproclaim
\bigskip

\demo{Proof}  We can find $h:\lambda^+ \rightarrow \lambda^+$ such that for
every $\gamma < \lambda^+$ the set \newline
$S_\gamma =: \{ \delta < \lambda^+:\text{cf}(\delta) = \lambda
\text{ and } h(\delta) = \gamma\}$ is stationary, so $\langle
S_\gamma:\gamma < \lambda \rangle$ is a partition \newline
of
$S =: \{ \delta < \lambda^+:\text{cf}(\delta) = \lambda\}$.  We can find
$\bar C^\gamma = \langle C_\delta:\delta \in S_\gamma \rangle$ such
that $C_\delta$ is a club of $\delta$ of order type $\lambda$.  For any
$\nu \in {}^{\omega >}(\lambda^+)$ we define:
\medskip
\roster
\item "{(a)}"  for $\ell < \ell g(\nu)$, if $\nu(\ell) \in S$ then let
\newline
$a_\ell = a_{\nu,\ell} = 
\{ \text{otp}(C_{\nu(\ell)} \cap \nu(k)):k < \ell g(\nu)
\text{ and } \nu(k) < \nu(\ell)\}$,
\item "{(b)}"  $\ell_\nu$ is the $\ell < \ell g(\nu)$ such that
{\roster
\itemitem{ (i) }  $\nu(\ell) \in S$,
\itemitem{ (ii) }  among those with $\sup(a_{\nu,\ell})$ is maximal, and
\itemitem{ (iii) }  among those with $\ell$ minimal,
\endroster}
\item "{(c)}"  if $\ell_\nu$ is well defined let $d(\nu) = h(\nu(\ell_\nu))$
otherwise let $d(\nu) = 0$.
\endroster
\medskip

Now suppose $\langle (u_\alpha,v_\alpha):\alpha < \lambda^+ \rangle,
\gamma < \lambda^+$ and $E$ are as in Definition \scite{2.1} 
and we shall prove the conclusion there.  Let \newline
$E^* = \{ \delta
\in E:\delta \text{ is a limit ordinal and } \alpha < \delta 
\Rightarrow \delta >$ \newline

$\qquad \qquad \qquad \sup [\cup \{ \text{Rang}(\eta):
\eta \in u_\alpha \cup v_\alpha\}]\}$. \newline
Clearly $E^* \subseteq E$ is a club of $\lambda^+$.
\medskip

For each $\delta \in S_\gamma$ let

$$
f_0(\delta) =: \sup[\delta \cap \bigcup \{ \text{Rang}(\nu):\nu \in 
u_\delta \cup v_\delta\}].
$$
\medskip

\noindent
As $\text{cf}(\delta) = \lambda > |u_\alpha \cup v_\alpha|$ and the
sequences are finite clearly $f_0(\delta) < \delta$.  Hence by Fodor's lemma
for some $\xi^*,S^1_\gamma =: \{ \delta \in S_\gamma:f_0(\delta) = \xi^*\}$
is a stationary subset of $\lambda^+$ (note that $\gamma$ is fixed here).  
Let $\xi^* = \dsize \bigcup_{i < \lambda}a_{2,i}$ where $a_{2,i}$ is 
increasing with $i$ and $|a_{2,i}| < \lambda$.  So for
$\delta \in S^1_\gamma$

$$
\align
f_1(\delta) = \text{ Min} \biggl\{ i < \lambda:&\delta \cap \bigcup
\{\text{Rang}(\nu):\nu \in u_\delta \cup v_\delta\} \\
  &\text{is a subset of } a_{2,i} \biggr\}
\endalign
$$
\medskip

\noindent
is a well defined ordinal $< \lambda$, hence for some $i^* < \lambda$
the set

$$
S^2_\gamma =: \{ \delta \in S^1_\gamma:f_1(\delta) = i^*\}
$$
\medskip

\noindent
is a stationary subset of $\lambda^+$.  For $\delta \in S^2_\gamma$ let

$$
\align
b_\delta =: \biggl\{ \text{otp}(C_\beta \cap \alpha):&\alpha < \beta \in
S \text{ and both} \\
  &\text{are in } a_{2,i^*} \cup \{ \delta \} \cup
\bigcup \{\text{Rang } \nu: \nu \in u_\delta \cup v_\delta\} \biggr\}.
\endalign
$$
\medskip

\noindent
So $b_\delta$ is a subset of $\lambda$ of cardinality $< \lambda$ hence
$\varepsilon_\delta =: \sup(b_\delta) < \lambda$, hence for some
$\varepsilon^*$

$$
S^3_\gamma =: \{ \delta \in S^2_\gamma:\varepsilon_\delta = \varepsilon^*\}
$$
\medskip

\noindent
is a stationary subset of $\lambda^+$.  Choose $\beta^*$ such that
\medskip
\roster
\item "{$(*)$}"  $\beta^* \in S^3_\gamma \cap E^* \text{ and }
\beta^* = \sup(\beta^* \cap S^3_\gamma \cap E^*)$.
\endroster
\medskip

\noindent
As $C_{\beta^*}$ has order type $\lambda$, (and is a club of $\beta^*$) for
some $\alpha^* \in \beta^* \cap S^3_\gamma \cap E^*$ we have
$\text{otp}(C_{\beta^*} \cap \alpha^*) > \varepsilon^*$. \newline
We want to show that $\alpha^*,\beta^*$ are as required.  Obviously
$\alpha^* < \beta^*,\alpha^* \in E$ and $\beta^* \in E$.  So assume
$\nu \in u_{\alpha^*},\rho \in v_{\beta^*}$ and we shall prove that
$d(\nu \char 94 \rho) = \gamma$, which suffices.  As $h(\beta^*) = \gamma$ (as
$\beta^* \in S^3_\gamma \subseteq S_\gamma)$ it suffices to prove that
$(\nu \char 94 \rho)(\ell_{\nu \char 94 \rho}) = \beta^*$.  Now for some
$\ell_0,\ell_1$ we have $\nu(\ell_0) = \alpha^*,\rho(\ell_1) = \beta^*$ (as
$\nu \in u_{\alpha^*},\rho \in v_{\beta^*})$ and since
$\text{otp}(C_{\beta^*} \cap
\alpha^*) > \varepsilon^*$, by the definition of $\ell_{\nu \char 94 \rho}$
it suffices to prove
\medskip
\roster
\item "{$\bigotimes$}"  \underbar{if} $\ell,k < \ell g(\nu \char 94 \rho),
(\nu \char 94 \rho)(\ell) \in S,(\nu \char 94 \rho)(k) < (\nu \char 94 \rho)
(\ell)$ \underbar{then}
{\roster
\itemitem{ (i) }  $\text{otp}[C_{(\nu \char 94 \rho)(\ell)} \cap
(\nu \char 94 \rho)(k)] \le \varepsilon^*$ or
\itemitem{ (ii) }  $(\nu \char 94 \rho)(\ell) = \beta^*$.
\endroster}
\endroster
\medskip

\noindent
Assume $\ell,k$ satisfy the assumption of $\otimes$ and we shall show its
conclusion.
\enddemo
\bigskip

\noindent
\underbar{Case 1}:  If $(\nu \char 94 \rho)(\ell)$ and 
$(\nu \char 94 \rho)(k)$ belong to

$$
a_{2,i^*} \cup \{ \beta^* \} \cup \bigcup \{\text{Rang}(\eta):
\eta \in u_{\beta^*} \cup v_{\beta^*}\}
$$
\medskip

\noindent
then clause $(i)$ holds because
\medskip
\roster
\item "{$(\alpha)$}"  $\text{otp}(C_{(\nu \char 94 \rho)(\ell)} \cap
(\nu \char 94 \rho)(k)) \in b_{\beta^*}$ (see the definition of
$b_{\beta^*}$) and
\item "{$(\beta)$}"  $\sup(b_{\beta^*}) = \varepsilon_{\beta^*}$ (see the
definition of $\varepsilon_{\beta^*}$) and
\item "{$(\gamma)$}"  $\varepsilon_{\beta^*} = \varepsilon^*$ (as
$\beta^* \in S^3_\gamma$ and see the choice of $\varepsilon^*$ and
$S^3_\gamma$).
\endroster
\bigskip

\noindent
\underbar{Case 2}:  If $(\nu \char 94 \rho)(\ell)$ and 
$(\nu \char 94 \rho)(k)$ belong to

$$
a_{2,i^*} \cup \bigcup \{\text{Rang}(\eta):\eta \in u_{\alpha^*} \cup
v_{\alpha^*}\}
$$
\medskip

\noindent
then the proof is similar to the proof of the previous case.
\bigskip

\noindent
\underbar{Case 3}:  No previous case. \newline

So $(\nu \char 94 \rho)(\ell)$ and $(\nu \char 94 \rho)(k)$ are not in
$a_{2,i^*}$, hence (as $\{\nu,\rho\} \subseteq (u_{\alpha^*} \cup
v_{\beta^*})$, and $\{ \alpha^*,\beta^*\} \subseteq S^2_\gamma \subseteq
S^1_\gamma$) 

$$
m \in \{\ell,k\} \and m < \ell g(\nu) \Rightarrow (\nu \char 94 \rho)(m)
= \nu(m) \ge \alpha^*,
$$

$$
m \in \{\ell,k\} \and m \ge \ell g(\nu) \Rightarrow (\nu \char 94 \rho)(m)
= \rho(m - \ell g(\nu)) \ge \beta^*.
$$
\medskip

\noindent
As $\beta^* \in E^*$ and $\beta^* > \alpha^*$ clearly 
$\sup(\text{Rang}(\nu)) < \beta^*$, but also \newline 
$(\nu \char 94 \rho)(k) < (\nu \char 94 \rho)(\ell)$
(see $\bigotimes$). \newline

Together necessarily $k < \ell g(\nu),\nu(k) \in [\alpha^*,\beta^*)$,
$\ell \in [\ell g(\nu),\ell g(\nu) + \ell g(\rho))$ and 
$\rho(\ell - \ell g(\nu)) \in [\beta^*,\lambda^+)$.  If 
$\rho(\ell) = \beta^*$ then clause $(ii)$ of the conclusion
holds.  Otherwise necessarily $\nu(\ell) > \beta^*$ hence

$$
\align
\text{otp}(C_{(\nu \char 94 \rho)(\ell)}) \cap (\nu \char 94 \rho)(k)) = 
&\text{ otp}(C_{\rho(\ell - \ell g(\nu))} \cap \nu(k)) \\
  &\le \text{ otp}(C_{\rho(\ell - \ell g(\nu))} \cap \beta^*) \le 
\sup(a_{\beta^*}) \le \varepsilon^*
\endalign
$$
\medskip

\noindent
so clause $(i)$ of $\otimes$ holds. \hfill$\square_{\scite{4.1}}$
\bigskip

\demo{\stag{4.2} Conclusion}  For $\lambda$ regular,
$Pr_1(\lambda^{+2},\lambda^{+2},\lambda^{+2},\lambda)$ holds.
\enddemo
\bigskip

\demo{Proof}  By \scite{4.1} and \scite{2.2}(1). \hfill$\square_{\scite{4.2}}$
\enddemo
\bigskip

\definition{\stag{4.3} Definition}  1) Let $Pr_6(\lambda,\theta,\sigma)$ 
means that for some
$\Xi$, an unbounded subset of $\{ \tau:\tau < \sigma,\tau$ is a cardinal
(finite or infinite)$\}$, there is a $d:{}^{\omega >}
(\lambda \times \Xi) \rightarrow
\omega$ such that if $\gamma < \theta$ and $\tau \in \Xi$ are given and
$\langle (u_\alpha,v_\alpha):\alpha < \lambda \rangle$ satisfies
\medskip
\roster
%
\item "{(i)}"  $u_\alpha \subseteq {}^{\omega >}(\lambda \times \Xi)
\backslash {}^{2 \ge}(\lambda \times \Xi)$,
\item "{(ii)}"  $v_\alpha \subseteq {}^{\omega >}(\lambda \times \Xi)
\backslash {}^{2 \ge}(\lambda \times \Xi)$,
\item "{(iii)}"  $|u_\alpha| = |v_\alpha| = \tau$,
\item "{(iv)}"  $\nu \in u_\beta \Rightarrow \nu(\ell g(\nu)-1) = \langle
\gamma,\tau \rangle$,
\item "{(v)}"  $\rho \in u_\alpha \Rightarrow \rho(0) = \langle \gamma,\tau
\rangle$,
\item "{(vi)}"  $\eta \in u_\alpha \cup v_\alpha \Rightarrow (\exists \ell)
(\eta(\ell) = \langle \alpha,\tau \rangle )$
\endroster
\medskip

\noindent
\underbar{then} for some $\alpha < \beta$ we have

$$
\nu \in u_\beta \and \rho \in v_\alpha \Rightarrow (\nu \char 94 \rho)[d(\nu
\char 94 \rho)] = \langle \gamma,\tau \rangle.
$$
\medskip

\noindent
2) Let $Pr_6(\lambda,\sigma)$ means $Pr_6(\lambda,\lambda,\sigma)$.
\enddefinition
\bigskip

\demo{\stag{4.4} Fact}  
$Pr_6(\lambda,\lambda,\theta,\sigma),\theta \ge \sigma$
implies $Pr_6(\lambda,\theta,\sigma)$.
\enddemo
\bigskip

\demo{Proof}  Let $c$ be a function from ${}^{\omega >}\lambda$ to $\theta$
exemplifying $Pr_6(\lambda,\lambda,\theta,\sigma)$.  Let $e$ be a one to one
function from $\theta \times \Xi$ onto $\theta$.

Now we define a function $d$ from ${}^{\omega >}(\lambda \times \Xi)$ to
$\omega$:

$$
d(\nu) = \text{ Min}\{\ell:c(\langle e(\nu(m)):m < \ell g(\nu) \rangle) =
e(\nu(\ell))\}. 
$$
\medskip
${}$ \hfill$\square_{\scite{4.4}}$
\enddemo
\bigskip

\proclaim{\stag{4.5} Claim}  If $Pr_6(\lambda^+,\sigma),\lambda$ 
regular and $\sigma
\le \lambda$ \underbar{then} $Pr_1(\lambda^{+2},\lambda^{+2},\lambda^{+2},
\sigma)$.
\endproclaim
\bigskip

\demo{Proof}  Like the proof of \scite{1.1}.
\enddemo
\bigskip

\remark{\stag{4.6} Remark}  As in \scite{4.1}, \scite{4.2} 
we can prove that if $\mu > \text{ cf}(\mu) + \sigma$ then \newline
$Pr^6(\mu^+,\mu^+,\mu^+,\sigma)$, hence
$Pr_1(\mu^{+2},\mu^{+2},\mu^{+2},\sigma)$, but this does not give new
information.
\endremark
\newpage

\shlhetal

\newpage

REFERENCES.  
\bibliographystyle{lit-plain}
\bibliography{lista,listb,listx}

\enddocument

\bye

%% file: mathdefs.tex
\expandafter\ifx\csname mathdefs.tex\endcsname\relax
  \expandafter\gdef\csname mathdefs.tex\endcsname{}
\else \message{Hey!  Apparently you were trying to
  \string\input{mathdefs.tex} twice.   This does not make sense.} 
\errmessage{Please edit your file (probably \jobname.tex) and remove
any duplicate ``\string\input'' lines} \fi




\catcode`\X=12\catcode`\@=11

\def\n@wcount{\alloc@0\count\countdef\insc@unt}
\def\n@wwrite{\alloc@7\write\chardef\sixt@@n}
\def\n@wread{\alloc@6\read\chardef\sixt@@n}
\def\r@s@t{\relax}\def\v@idline{\par}\def\@mputate#1/{#1}
\def\l@c@l#1X{\firstpart.#1}\def\gl@b@l#1X{#1}\def\t@d@l#1X{{}}

\def\crossrefs#1{\ifx\all#1\let\tr@ce=\all\else\def\tr@ce{#1,}\fi
   \n@wwrite\cit@tionsout\openout\cit@tionsout=\jobname.cit 
   \write\cit@tionsout{\tr@ce}\expandafter\setfl@gs\tr@ce,}
\def\setfl@gs#1,{\def\@{#1}\ifx\@\empty\let\next=\relax
   \else\let\next=\setfl@gs\expandafter\xdef
   \csname#1tr@cetrue\endcsname{}\fi\next}
\def\m@ketag#1#2{\expandafter\n@wcount\csname#2tagno\endcsname
     \csname#2tagno\endcsname=0\let\tail=\all\xdef\all{\tail#2,}
   \ifx#1\l@c@l\let\tail=\r@s@t\xdef\r@s@t{\csname#2tagno\endcsname=0\tail}\fi
   \expandafter\gdef\csname#2cite\endcsname##1{\expandafter
     \ifx\csname#2tag##1\endcsname\relax?\else\csname#2tag##1\endcsname\fi
     \expandafter\ifx\csname#2tr@cetrue\endcsname\relax\else
     \write\cit@tionsout{#2tag ##1 cited on page \folio.}\fi}
   \expandafter\gdef\csname#2page\endcsname##1{\expandafter
     \ifx\csname#2page##1\endcsname\relax?\else\csname#2page##1\endcsname\fi
     \expandafter\ifx\csname#2tr@cetrue\endcsname\relax\else
     \write\cit@tionsout{#2tag ##1 cited on page \folio.}\fi}
   \expandafter\gdef\csname#2tag\endcsname##1{\expandafter
      \ifx\csname#2check##1\endcsname\relax
      \expandafter\xdef\csname#2check##1\endcsname{}%
      \else\immediate\write16{Warning: #2tag ##1 used more than once.}\fi
      \multit@g{#1}{#2}##1/X%
      \write\t@gsout{#2tag ##1 assigned number \csname#2tag##1\endcsname\space
      on page \number\count0.}%
   \csname#2tag##1\endcsname}}
\def\multit@g#1#2#3/#4X{\def\t@mp{#4}\ifx\t@mp\empty%
      \global\advance\csname#2tagno\endcsname by 1 
      \expandafter\xdef\csname#2tag#3\endcsname
      {#1\number\csname#2tagno\endcsnameX}%
   \else\expandafter\ifx\csname#2last#3\endcsname\relax
      \expandafter\n@wcount\csname#2last#3\endcsname
      \global\advance\csname#2tagno\endcsname by 1 
      \expandafter\xdef\csname#2tag#3\endcsname
      {#1\number\csname#2tagno\endcsnameX}
      \write\t@gsout{#2tag #3 assigned number \csname#2tag#3\endcsname\space
      on page \number\count0.}\fi
   \global\advance\csname#2last#3\endcsname by 1
   \def\t@mp{\expandafter\xdef\csname#2tag#3/}%
   \expandafter\t@mp\@mputate#4\endcsname
   {\csname#2tag#3\endcsname\lastpart{\csname#2last#3\endcsname}}\fi}
\def\t@gs#1{\def\all{}\m@ketag#1e\m@ketag#1s\m@ketag\t@d@l p
   \m@ketag\gl@b@l r \n@wread\t@gsin
   \openin\t@gsin=\jobname.tgs \re@der \closein\t@gsin
   \n@wwrite\t@gsout\openout\t@gsout=\jobname.tgs }
\outer\def\localtags{\t@gs\l@c@l}
\outer\def\globaltags{\t@gs\gl@b@l}
\outer\def\newlocaltag#1{\m@ketag\l@c@l{#1}}
\outer\def\newglobaltag#1{\m@ketag\gl@b@l{#1}}

\newif\ifpr@ 
\def\m@kecs #1tag #2 assigned number #3 on page #4.%
   {\expandafter\gdef\csname#1tag#2\endcsname{#3}
   \expandafter\gdef\csname#1page#2\endcsname{#4}
   \ifpr@\expandafter\xdef\csname#1check#2\endcsname{}\fi}
\def\re@der{\ifeof\t@gsin\let\next=\relax\else
   \read\t@gsin to\t@gline\ifx\t@gline\v@idline\else
   \expandafter\m@kecs \t@gline\fi\let \next=\re@der\fi\next}
\def\pretags#1{\pr@true\pret@gs#1,,}
\def\pret@gs#1,{\def\@{#1}\ifx\@\empty\let\n@xtfile=\relax
   \else\let\n@xtfile=\pret@gs \openin\t@gsin=#1.tgs \message{#1} \re@der 
   \closein\t@gsin\fi \n@xtfile}

\newcount\sectno\sectno=0\newcount\subsectno\subsectno=0
\newif\ifultr@local \def\ultralocal{\ultr@localtrue}
\def\firstpart{\number\sectno}
\def\lastpart#1{\ifcase#1 \or a\or b\or c\or d\or e\or f\or g\or h\or 
   i\or k\or l\or m\or n\or o\or p\or q\or r\or s\or t\or u\or v\or w\or 
   x\or y\or z \fi}

\def\resetall{\global\advance\sectno by 1\subsectno=0
   \gdef\firstpart{\number\sectno}\r@s@t}
\def\resetsub{\global\advance\subsectno by 1
   \gdef\firstpart{\number\sectno.\number\subsectno}\r@s@t}
\def\newsection#1\par{\resetall\vskip0pt plus.3\vsize\penalty-250
   \vskip0pt plus-.3\vsize\bigskip\bigskip
   \message{#1}\leftline{\bf#1}\nobreak\bigskip}
\def\subsection#1\par{\ifultr@local\resetsub\fi
   \vskip0pt plus.2\vsize\penalty-250\vskip0pt plus-.2\vsize
   \bigskip\smallskip\message{#1}\leftline{\bf#1}\nobreak\medskip}

\def\t@gsoff#1,{\def\@{#1}\ifx\@\empty\let\next=\relax\else\let\next=\t@gsoff
   \def\@@{p}\ifx\@\@@\else
   \expandafter\gdef\csname#1cite\endcsname##1{\zeigen{##1}}
   \expandafter\gdef\csname#1page\endcsname##1{?}
   \expandafter\gdef\csname#1tag\endcsname##1{\zeigen{##1}}\fi\fi\next}
\def\verbatimtags{\ifx\all\relax\else\expandafter\t@gsoff\all,\fi}
\def\zeigen#1{\hbox{$\langle$}#1\hbox{$\rangle$}}

\def\(#1){\edef\dot@g{\ifmmode\ifinner(\hbox{\noexpand\etag{#1}})
   \else\noexpand\eqno(\hbox{\noexpand\etag{#1}})\fi
   \else(\noexpand\ecite{#1})\fi}\dot@g}

\newif\ifbr@ck
\def\eat#1{}
\def\[#1]{\br@cktrue[\br@cket#1'X]}
\def\br@cket#1'#2X{\def\temp{#2}\ifx\temp\empty\let\next\eat
   \else\let\next\br@cket\fi
   \ifbr@ck\br@ckfalse\br@ck@t#1,X\else\br@cktrue#1\fi\next#2X}
\def\br@ck@t#1,#2X{\def\temp{#2}\ifx\temp\empty\let\neext\eat
   \else\let\neext\br@ck@t\def\temp{,}\fi
   \def\teemp{#1}\ifx\teemp\empty\else\rcite{#1}\fi\temp\neext#2X}
\def\resetbr@cket{\gdef\[##1]{[\rtag{##1}]}}
\def\references{\resetbr@cket\newsection References\par}

\newtoks\symb@ls\newtoks\s@mb@ls\newtoks\p@gelist\n@wcount\ftn@mber
    \ftn@mber=1\newif\ifftn@mbers\ftn@mbersfalse\newif\ifbyp@ge\byp@gefalse
\def\defm@rk{\ifftn@mbers\n@mberm@rk\else\symb@lm@rk\fi}
\def\n@mberm@rk{\xdef\m@rk{{\the\ftn@mber}}%
    \global\advance\ftn@mber by 1 }
\def\rot@te#1{\let\temp=#1\global#1=\expandafter\r@t@te\the\temp,X}
\def\r@t@te#1,#2X{{#2#1}\xdef\m@rk{{#1}}}
\def\b@@st#1{{$^{#1}$}}\def\str@p#1{#1}
\def\symb@lm@rk{\ifbyp@ge\rot@te\p@gelist\ifnum\expandafter\str@p\m@rk=1 
    \s@mb@ls=\symb@ls\fi\write\f@nsout{\number\count0}\fi \rot@te\s@mb@ls}
\def\byp@ge{\byp@getrue\n@wwrite\f@nsin\openin\f@nsin=\jobname.fns 
    \n@wcount\currentp@ge\currentp@ge=0\p@gelist={0}
    \re@dfns\closein\f@nsin\rot@te\p@gelist
    \n@wread\f@nsout\openout\f@nsout=\jobname.fns }
\def\m@kelist#1X#2{{#1,#2}}
\def\re@dfns{\ifeof\f@nsin\let\next=\relax\else\read\f@nsin to \f@nline
    \ifx\f@nline\v@idline\else\let\t@mplist=\p@gelist
    \ifnum\currentp@ge=\f@nline
    \global\p@gelist=\expandafter\m@kelist\the\t@mplistX0
    \else\currentp@ge=\f@nline
    \global\p@gelist=\expandafter\m@kelist\the\t@mplistX1\fi\fi
    \let\next=\re@dfns\fi\next}
\def\symbols#1{\symb@ls={#1}\s@mb@ls=\symb@ls} 
\def\bigsymbol{\textstyle}
\symbols{\bigsymbol\ast,\dagger,\ddagger,\sharp,\flat,\natural,\star}
\def\ftnumbers{\ftn@mberstrue} \def\ftsymbols{\ftn@mbersfalse}
\def\paginal{\byp@ge} \def\resetftnumbers{\ftn@mber=1}
\def\ftnote#1{\defm@rk\expandafter\expandafter\expandafter\footnote
    \expandafter\b@@st\m@rk{#1}}

\long\def\jump#1\endjump{}
\def\ssum{\mathop{\lower .1em\hbox{$\textstyle\Sigma$}}\nolimits}

\def\qed{\nobreak\kern 1em \vrule height .5em width .5em depth 0em}
\def\newneq{\hbox{\rlap{\hbox to 1\wd9{\hss$=$\hss}}\raise .1em 
   \hbox to 1\wd9{\hss$\scriptscriptstyle/$\hss}}}
\def\subsetne{\setbox9 = \hbox{$\subset$}\mathrel{\hbox{\rlap
   {\lower .4em \newneq}\raise .13em \hbox{$\subset$}}}}
\def\supsetne{\setbox9 = \hbox{$\subset$}\mathrel{\hbox{\rlap
   {\lower .4em \newneq}\raise .13em \hbox{$\supset$}}}}

\def\vbar{\mathchoice{\vrule height6.3ptdepth-.5ptwidth.8pt\kern-.8pt}
   {\vrule height6.3ptdepth-.5ptwidth.8pt\kern-.8pt}
   {\vrule height4.1ptdepth-.35ptwidth.6pt\kern-.6pt}
   {\vrule height3.1ptdepth-.25ptwidth.5pt\kern-.5pt}}
\def\f@dge{\mathchoice{}{}{\mkern.5mu}{\mkern.8mu}}
\def\b@c#1#2{{\rm \mkern#2mu\vbar\mkern-#2mu#1}}
\def\b@b#1{{\rm I\mkern-3.5mu #1}}
\def\b@a#1#2{{\rm #1\mkern-#2mu\f@dge #1}}
\def\bb#1{{\count4=`#1 \advance\count4by-64 \ifcase\count4\or\b@a A{11.5}\or
   \b@b B\or\b@c C{5}\or\b@b D\or\b@b E\or\b@b F \or\b@c G{5}\or\b@b H\or
   \b@b I\or\b@c J{3}\or\b@b K\or\b@b L \or\b@b M\or\b@b N\or\b@c O{5} \or
   \b@b P\or\b@c Q{5}\or\b@b R\or\b@a S{8}\or\b@a T{10.5}\or\b@c U{5}\or
   \b@a V{12}\or\b@a W{16.5}\or\b@a X{11}\or\b@a Y{11.7}\or\b@a Z{7.5}\fi}}

\catcode`\X=11 \catcode`\@=12